\documentclass[11pt, a4paper,reqno]{amsart}
\usepackage{preemble-art}

\begin{document}

\title[Topological flows]
      {Topological flows for hyperbolic groups}
\author{Ryokichi Tanaka}
\address{Mathematical  Institute, Tohoku  University, Sendai 980-8578 JAPAN}
\email{ryokichi.tanaka@tohoku.ac.jp}
\date{\today}

\maketitle

\begin{abstract}
We show that for every non-elementary hyperbolic group the Bowen-Margulis current associated with a strongly hyperbolic metric forms a unique group-invariant Radon measure class of maximal Hausdorff dimension on the boundary square.
Applications include a characterization of roughly similar hyperbolic metrics via mean distortion. 
\end{abstract}

\section{Introduction}

Let $\G$ be a non-elementary hyperbolic group.
It is known that the ideal boundary $\partial \G$ admits plenty of finite measures which give rise to $\G$-invariant measures on the boundary square $\partial^2 \G:=(\partial \G)^2\setminus \{{\rm diagonal}\}$, where
we consider the diagonal action of $\G$ on $\partial^2 \G$.
The $\G$-invariant measures on $\partial^2 \G$ arise in the study of geodesic flows on negatively curved manifolds since there exists a natural correspondence between geodesic flow-invariant finite measures and $\G$-invariant Radon (i.e., locally finite and Borel regular) measures on the boundary square of the universal coverings. 
For example, one way to obtain an invariant measure of maximal entropy for the geodesic flow on a compact negatively curved manifold is to construct a $\G$-invariant Radon measure in the measure class consisting of the product of Patterson-Sullivan measures on $\partial^2 \G$.
This construction has been generalized beyond the manifold setting such as $\CAT(-1)$-spaces \cite{BourdonStructure}.
For general hyperbolic groups,
the study of geodesic flows goes back to the original paper by Gromov \cite[8.3]{GromovHyperbolic}.
Furman \cite{FurmanCoarse} has given a framework for a general hyperbolic group $\G$ to analyze $\G$-invariant measures on $\partial^2 \G$ without assuming any actions on another reasonable geometric space (see also \cite{BaderFurman}).
It would be desirable to define the entropy of a $\G$-invariant measure on $\partial^2 \G$ for a general hyperbolic group such that it gives the measure theoretical entropy for a geodesic flow when the geodesic flow is defined.
Instead we discuss the Hausdorff dimension; the entropy coincides with the Hausdorff dimension under an appropriate normalization
in the case of negatively curved manifolds as it was studied and the approach has been suggested for hyperbolic groups by Kaimanovich \cite[Section 3.5]{KaimanovichInvariant}.
The purpose of this paper is to characterize the $\G$-invariant Radon measure class where the maximal entropy is replaced by the maximal Hausdorff dimension.

For any non-elementary hyperbolic group $\G$ (not necessarily torsion-free), we consider a strongly hyperbolic metric $\wh d$ which is left-invariant and quasi-isometric to a word metric on $\G$.
For example, one may take the hat metric constructed by Mineyev \cite{MineyevFlow} and any hyperbolic group admits such a metric.
Let us consider the Bowen-Margulis current for $\wh d$,
\[
\wh \L:=\exp(2 \wh v\, \langle \x | \y\rangle_o)\,\wh \m_o\otimes \wh \m_o \quad \text{on $\partial^2 \G$},
\]
where $\langle \x|\y\rangle_o$ is the Gromov product based at $o$ and $\wh v$ is the exponential volume growth rate of $(\G, \wh d\,)$ and $\wh \m_o$ is the corresponding Patterson-Sullivan measure on $\partial \G$.
Letting the gauge $\wh \rho(\x, \y):=\exp(-\langle \x|\y \rangle_o)$ and $\wh \rho_\times$ be the maximum of $\wh \rho$ on each component on $\partial^2\G$,
we define the (lower) Hausdorff dimension of a Borel measure $\L$ on $\partial^2 \G$ by
\[
\underline \dim_H(\L, \wh \rho_\times):=\inf\{\dim_H (A, \wh \rho_\times) \ : \ \L(A) >0\}.
\]
It is known that $\dim_H(\partial^2 \G, \wh \rho_\times)=2\wh v$.
We show that the Bowen-Margulis current $\wh \L$ is the unique $\G$-invariant Radon measure on $\partial^2 \G$ of maximal Hausdorff dimension up to a multiplicative constant.

\begin{theorem}\label{Thm:uniqueMMD_hat}
Let $\G$ be a non-elementary hyperbolic group and $\wh d$ be a left-invariant strongly hyperbolic metric quasi-isometric to a word metric.
For any $\G$-invariant Radon measure $\Lambda$ on $\partial^2 \G$, 
if $\underline{\dim}_H(\Lambda, \wh \rho_\times)=\dim_H(\partial^2 \G, \wh \rho_\times)$, 
then $\Lambda$ is a constant multiple of the Bowen-Margulis current $\wh \L$.
\end{theorem}

The result allows us to compare any $\G$-invariant Radon measure $\L$ with the Bowen-Margulis current on $\partial^2 \G$.
Our proof is based on a topological model of geodesic flow on a general hyperbolic space introduced by Mineyev \cite{MineyevFlow}.
Note that one is able to compare two Bowen-Margulis currents arising from two hyperbolic metrics on $\G$ in a measure theoretical model introduced by Furman \cite{FurmanCoarse} and by the double ergodicity result of Kaimanovich \cite[Theorem 2.8]{KaimanovichErgodicity}, Bader and Furman \cite{BaderFurman} (see also \cite{Garncarek}) and Coulon {\it et al}.\ \cite[Section 4]{CDST}.
However, since our aim is to deal with any $\G$-invariant Radon measure $\L$ and the question involves the Hausdorff dimension, 
we require a topological setting independent of the choice of measures.

In the case of free groups with actions on metric trees, Kapovich and Nagnibeda have introduced the {\it geometric entropy} of a current
\cite[Definition 3.1]{KapovichNagnibedaGeometric}.
Analogously, we define
\[
h(\L, \wh d\,):=\lim_{R \to \infty}\liminf_{\wh d(x, y) \to \infty \ x, y \in \G}-\frac{\log \L \(O_o(x, R) \times O_o(y, R)\)}{\wh d(x, y)},
\]
where $O_o(x, R)$ is the shadow based at $o$ of the ball centered at $x$ with radius $R$ (Section \ref{Sec:shadows})
and the limit exists as $R \to \infty$ since it is non-decreasing.
The lower Hausdorff dimension of $\L$ is given by for a large enough $R>0$,
\[
\underline \dim_H (\L, \wh \rho_\times)=\inf_{\text{$\L$-a.e.\ $(\x, \y)$}}\liminf_{n \to \infty}-\frac{\log\L \(O_o(\phi(-n), R) \times O_o(\phi(n), R)\)}{n},
\]
where $\phi$ is a rough geodesic with extreme points $(\x, \y)$ normalized so that $\phi(0)$ realizes a minimum distance to the base point $o$ and 
$\inf_{\text{$\L$-a.e.\ $(\x, \y)$}}$ stands for the essential infimum relative to $\L$ (Section \ref{Sec:dimension}).
Then we have 
\[
2h(\L, \wh d\,)\le \underline \dim_H (\L, \wh \rho_\times).
\]
In the case of (simplicial) trees, their result reads if $h(\L, \wh d\,)=\wh v$, then $\L$ is a constant multiple of the Bowen-Margulis current (in our terminology) \cite[Theorems C and E]{KapovichNagnibedaGeometric}.
Since geodesic metrics in trees are strongly hyperbolic (Section \ref{Sec:strongly_hyperbolic}), Theorem \ref{Thm:uniqueMMD_hat} reproduces the same result.
On the other hand, we note that one is not able to replace $\underline \dim_H(\L, \wh d\,)$ by the upper Hausdorff dimension $\overline \dim_H(\L, \wh d\,)$ (Section \ref{Sec:dimension}) to obtain Theorem \ref{Thm:uniqueMMD_hat} unless $\L$ is ergodic with respect to the $\G$-action, i.e., 
any $\G$-invariant Borel set $A$ in $\partial^2 \G$, one has $\L(A)=0$ or $\L(\partial^2 \G \setminus A)=0$.

\subsection{A connection to mean distortion}
We have the corresponding result for word metrics (Theorem \ref{Thm:uniqueMMD_word}).
Let us illustrate an application to comparison between two word metrics.
For any two finite symmetric sets of generators $S$, $S^\star$, we consider the word norms $|\cdot|_S$ and $|\cdot|_{S^\star}$ for $S$ and $S^\star$, respectively.
For each integer $n \ge 0$, let 
\[
\bS_n:=\{x \in \G \ : \ |x|_S=n\}.
\]
We sample an element $x_n$ according to the uniform distribution $\Unif_{S, n}$ on $\bS_n$ for each $n$, and consider the (liminf-) linear growth rate of the average length with respect to 
$|\cdot|_{S^\star}$,
\[
\t(S^\star/ S):=\liminf_{n \to \infty}\frac{1}{n}\E_{\Unif_{S, n}}|x_n|_{S^\star},  \quad \text{where} \quad \E_{\Unif_{S, n}}|x_n|_{S^\star}:=\frac{1}{|\bS_n|}\sum_{x \in \bS_n}|x|_{S^\star}.
\]
Let us call $\t(S^\star/S)$ the {\it mean distortion} of $|\cdot|_{S^\star}$ relative to $|\cdot|_S$.
Since two word metrics are bi-Lipschitz, the mean distortion is bounded from above and from below by the bi-Lipschitz constants.
The quantity $\t(S^\star/S)$ represents a typical distortion rate between two word norms. 
If we denote by $\gr(S)$ and $\gr(S^\star)$ the exponential volume growth rates for word metrics relative to $S$ and $S^\star$, respectively,
then a simple counting argument yields
$\t(S^\star/S) \ge \gr(S)/\gr(S^\star)$.
We show that the equality holds if and only if two word metrics are roughly similar.

\begin{theorem}\label{Thm:mean_distortion}
Let $\G$ be a non-elementary hyperbolic group.
For any pair of finite symmetric sets of generators $S$ and $S^\star$ in $\G$,
we have that
\[
\t(S^\star/ S)=\lim_{n \to \infty}\frac{1}{n}\E_{\Unif_{S, n}}|x_n|_{S^\star},
\]
and
for any $\e>0$, 
\begin{equation*}
\lim_{n \to \infty} \Unif_{S, n}\left\{x \in {\bf S}_n \ : \ \frac{|d_{S^\star}(o, x)-n \t(S^\star/S)|}{n} > \e\right\}=0.
\end{equation*}
Moreover,
it holds that
\[
\t(S^\star/S) \ge \frac{\gr(S)}{\gr(S^\star)},
\]
and the equality $\t(S^\star/S)=\gr(S)/\gr(S^\star)$ holds if and only if word metrics $d_{S^\star}$ and $d_S$ are roughly similar, i.e.,
there exist constants $\t>0$ and $D \ge 0$ such that
\[
|d_{S^\star}(x, y)-\t d_S(x, y)| \le D \quad \text{for all $x, y \in \G$}.
\]
\end{theorem}

Theorem \ref{Thm:mean_distortion} is a consequence of Theorems \ref{Thm:distortion} and \ref{Thm:distortion-similarity}.
We have stated Theorem \ref{Thm:mean_distortion} only for word metrics; but this is mainly for the sake of simplicity on the statement---one may discuss more general hyperbolic metrics.
The proof is based on showing that the Hausdorff dimension computed by the gauge associated with $d_{S^\star}$ of the Bowen-Margulis current for $d_{S}$ coincides with 
$\gr(S)/\t(S^\star/S)$.
The inequality $\t(S^\star/S) \ge \gr(S)/\gr(S^\star)$ has been obtained by Calegari and Fujiwara \cite[Remark 4.28]{CalegariFujiwara2010}
(and this remark has motivated our result).
Moreover, $\t(S^\ast/S)$ is an algebraic number \cite[Corollary 4.27]{CalegariFujiwara2010}.
In fact, they have shown a central limit theorem which implies that $|x_n|_{S^\star}=\t(S^\star/S)n+O(\sqrt{n})$ with high probability for uniformly chosen $x_n$ in the sphere $\bS_n$ as $n \to \infty$ (see also \cite[Corollary 3.6.4]{Calegari}).
Our proof indicates that the $O(\sqrt{n})$-fluctuation can not be negligible as soon as two word metrics are not roughly similar.
In the case of free groups and word metrics associated with free bases, the characterization of rough similarity has been given by \cite[Theorem F]{KKS} (where the mean distortion is called the generic stretching factor).
We point out that Theorem \ref{Thm:mean_distortion} is regarded as a discrete counterpart of a result for compact negatively curved manifolds by Knieper \cite[Theorem 1.2]{KnieperVolume} (where $\t(S^\star/S)$ corresponds to the geodesic stretch).

\subsection{Outline of the proof}
Let us give an outline of the proof of Theorem \ref{Thm:uniqueMMD_hat}.
We consider the flow $\{\wt \F_t\}_{t \in \R}$ on the space $\partial^2 \G \times \R$
defined by the translation in the $\R$-coordinate.
The space $\partial^2 \G \times \R$ admits a $\G$-action, which is constructed as in the following.
For a strongly hyperbolic metric $\wh d$ on $\G$, the corresponding Gromov product and the Busemann functions $\wh b_o(x, \x)$ extend continuously on the compactification $\G \cup \partial \G$.
By using the cocycle 
\[
\a(x, \x, \y):=\frac{1}{2}(\wh b_o(x^{-1}, \x)-\wh b_o(x^{-1}, \y)) \quad \text{for $x \in \G, \ (\x, \y) \in \partial^2 \G$},
\]
we define the $\G$-action on $\partial^2 \G \times \R$ by $x\cdot(\x, \y, t):=(x \x, x \y, t-\a(x, \x, \y))$.
Let us call this action the $(\G, \a)$-action.
The $(\G, \a)$-action is properly discontinuous and cocompact (Lemma \ref{Lem:Gamma-action}).
Then we define $\Fc_\a:=\G\backslash (\partial^2 \G \times \R)$ and call it the {\it topological flow space}.
Since the flow $\{\wt \F_t\}_{t \in \R}$ commutes with the $(\G, \a)$-action, it descends to the flow $\{\F_t\}_{t \in \R}$ on $\Fc_\a$.
It amounts to consider the geodesic flow on the total space of a unit tangent bundle in the case of manifolds.

This flow $\{\F_t\}_{t \in \R}$ on the topological flow space $\Fc_\a$ has potentially a lot of similarities to Axiom $A$ flows.
However, a direct connection to the known machinery seems lacking (e.g., \cite[Section 6.3]{BridgemanCanarySambarino}).
Moreover, the $(\G, \a)$-action on $\partial^2 \G \times \R$ is not necessarily free, and $\Fc_\a$ is far from being a manifold, 
it is not clear that one could resemble the techniques developed in the manifold setting;
in particular, a serious issue arises when one tries to obtain a lower bound of the topological entropy of the flow on $\Fc_\a$.
Instead, we follow a classical approach in the Axiom $A$ flows by Bowen and Ruelle \cite{BowenRuelleAxiomAflows}.
We construct a symbolic coding of the topological flow space $\Fc_\a$ by using Cannon's automatic structure of hyperbolic group $\G$.
This allows us to work with a two-sided subshift of finite type $(\SS, \s)$ based on the underlying graph structure.
Then the suspension flow $\Sus(\SS, r)$ with a natural roof function associated with the cocycle $\a$ carries a coding map
$w_\ast:\Sus(\SS, r) \to \Fc_\a$, which is bounded-to-one and equivariant with flows.
The suspension flow $\Sus(\SS, r)$ is now the place to work on, but the problem is that it is not clear as to whether the coding map is one-to-one over a residual set, and the two-sided shift space $(\SS, \s)$ is not necessarily topologically transitive.
In the $\CAT(-1)$-setting, a coding as satisfactory as in hyperbolic basic sets has been constructed in \cite{CLTstrong},
but we do not know as to whether a similar construction is possible for a general hyperbolic group.
An advantage, however, to use Cannon's automatic structure for coding is that it respects the geometry of a Cayley graph.
It is actually strong enough to employ thermodynamic formalism on $(\SS, \s)$.
The Bowen-Margulis current for a strongly hyperbolic metric is encoded on the shift space by using a H\"older continuous potential associated with the Busemann cocycle.
It requires a careful treatment to understand the support of this encoded measure since $(\SS, \s)$ is not transitive and we use a spectral decomposition into transitive components.
We show a key proposition which makes it possible to encode all $\G$-invariant Radon measures on $\partial^2 \G$ as a shift-invariant probability measure on $\SS$.
Then the problem is basically reduced to the uniqueness of measure of maximal entropy on a transitive subshift of finite type.

Let us mention a direction which we have not pursued in this paper;
given all this framework, it would be interesting to extend results to the product of Patterson-Sullivan measures associated with a H\"older continuous cocycle (e.g.,
\cite{LedrappierStructure}) for a general hyperbolic group setting.

\subsection{Organization of the paper}
In Section \ref{Sec:preliminary}, we review basics on hyperbolic metrics and Patterson-Sullivan measures.
In particular, we discuss a strongly hyperbolic metric and an associated cocycle we use.
We also give basic facts on the Hausdorff dimension of sets and measures on the boundary of a hyperbolic group.
In Section \ref{Sec:topological_flow_space}, we construct a topological flow space $\Fc_\a$ and show that the Bowen-Margulis current yields an ergodic flow-invariant measure on $\Fc_\a$ (Theorem \ref{Thm:ergodic}).
This implies the double ergodicity of the Patterson-Sullivan measures (Corollary \ref{Cor:double-ergodicity}).
We follow \cite{BaderFurman} and \cite[Appendix A]{Garncarek} (see also Coulon {\it et al}.\ \cite[Section 4]{CDST}) for the proofs in our setting.
In Section \ref{Sec:coding}, we construct a two-sided subshift of finite type based on an automatic structure of the hyperbolic group.
We show that there exists a natural coding map from the suspension flow $\Sus(\SS, r)$ to $\Fc_\a$ (Proposition \ref{Prop:coding}) and a key proposition which states that all $\G$-invariant Radon measures on $\partial^2 \G$ arise from shift-invariant probability measures on $\SS$ (Proposition \ref{Prop:coding_measure}).
In Section \ref{Sec:symbolic}, we use thermodynamic formalism to construct a shift-invariant measure which induces a measure dominated by the Bowen-Margulis current up to a positive multiplicative constant (Lemma \ref{Lem:coding_PS}).
We also formulate the variational principle based on the subshift of finite type $(\SS, \s)$.
In Section \ref{Sec:entropy}, we give a local entropy-dimension estimate which gives a direct connection between the measure theoretic entropy of an invariant probability measure on the shift space $(\SS, \s)$ and the Hausdorff dimension of a measure on $\partial^2 \G$ in the case when the two measures are connected as in Proposition \ref{Prop:coding_measure} (Lemma \ref{Lem:Local_dimension}).
We prove Theorem \ref{Thm:uniqueMMD_main} (i.e., Theorem \ref{Thm:uniqueMMD_hat}) for a strongly hyperbolic metric.
In Section \ref{Sec:mean_distortion}, we show the corresponding result for a word metric (Theorem \ref{Thm:uniqueMMD_word}).
We prove results on mean distortions (Theorems \ref{Thm:distortion} and \ref{Thm:distortion-similarity}) and deduce Theorem \ref{Thm:mean_distortion}.

{\bf Notation}:
Throughout this article, we write numerical constants $C, C', C'', \dots$ whose exact values may change from lines to lines, and
we denote by $C_\d$, etc.\ to indicate its dependance on $\d$ for a parameter $\d$.
For two real-valued functions $f(t)$ and $g(t)$,
we write $f(t)=g(t)+O_\d$ if and only if there exists a constant $C_\d$ such that $|f(t)-g(t)| \le C_\d$ for all $t \in I$.
Equivalently, we also sometimes write $f(t)=g(t)\pm C_\d$ if we want to emphasize the dependence on constants slightly more.
Although we might not state it everywhere, all measures on topological spaces are Borel.

\section{Preliminary}\label{Sec:preliminary}

\subsection{Metrics on hyperbolic groups}
\subsubsection{Hyperbolic groups}
Let $\G$ be a finitely generated group.
We say that a metric $d$ on $\G$ is {\it left-invariant} if $d(\g x, \g y)=d(x, y)$ for all $\g, x, y \in \G$,
and $d$ is $\d$-{\it hyperbolic} if there exists $\d \ge 0$ such that for any $x, y, z, w \in \G$,
\[
(x|y)_w \ge \min\{(x|z)_w, (z|y)_w\}-\d,
\]
where we define the Gromov product
\[
(x|y)_w:=\frac{d(w, x)+d(w, y)-d(x, y)}{2}.
\]
A metric $d$ is called {\it hyperbolic} if $d$ is $\d$-hyperbolic for some $\d \ge 0$.
For a finite set of generators $S$ with $S=S^{-1}$ in $\G$, the associated word metric is defined by for $x, y \in \G$,
\[
d_S(x, y):=|x^{-1}y|_S \quad \text{where $|x|_S:=\min\{k \ge 0 \ : \ s_1\cdots s_k=x, \ s_1, \dots, s_k \in S\}$}
\]
and $|\id|_S=0$ where $\id$ denotes the identity element.
We denote the associated Cayley graph by $\Cay(\G, S)$ and regard $d_S$ as the graph distance in $\Cay(\G, S)$.
The word metric $d_S$ is left-invariant, and it is proper, i.e., every ball of finite radius consists of finitely many points.
A finitely generated group $\G$ is called a {\it hyperbolic group} (or a {\it word hyperbolic group}) if a word metric is $\d$-hyperbolic for some $\d \ge 0$.
In fact, if $\G$ is a hyperbolic group, then any word metric is hyperbolic although the constant $\d$ depends on the set of generators $S$.
This follows since for every word metric any two points are joined by an isometric image of a (discrete) path, and any two word metrics are quasi-isometric.

Recall that two metrics $d_1$ and $d_2$ on $\G$ are {\it quasi-isometric} if for some constants $L>0, C \ge 0$ we have that
\[
(1/L)d_1(x, y) -C \le d_2(x, y) \le L d_1(x, y)+C \quad \text{for all $x, y \in \G$}.
\]
Let $I$ be an interval in $\R$.
We say that a path $\phi:I \to (\G, d)$ is an $(L, C)$-{\it quasi-geodesic} for constants $L, C$
if 
\[
(1/L)|s-t|-C \le d(\phi(s), \phi(t)) \le L|s-t|+C \quad \text{for all $s, t \in I$},
\] 
and $\phi:I \to (\G, d)$ is a $C$-{\it rough geodesic} 
if 
\[
|s-t| -C \le d(\phi(s), \phi(t)) \le |s-t| +C \quad \text{for all $s, t \in I$},
\]
where we say that $\phi$ is {\it geodesic} if $d(\phi(s), \phi(t))=|s-t|$ for all $s, t \in I$.
Note that if $d_1$ and $d_2$ are quasi-isometric,
then
geodesic or $C$-rough geodesic paths into $(\G, d_1)$ are $(L, C')$-quasi-geodesics in $(\G, d_2)$ for some $L, C'$.
A metric space $(\G, d)$ is called {\it geodesic} if any two points are joined by a geodesic path, and $(\G, d)$ is called $C$-{\it roughly geodesic} if for any two points $x, y \in \G$ there exists a $C$-rough geodesic path $\phi:[a, b] \to (\G, d)$ such that $\phi(a)=x$ and $\phi(b)=y$.
We say that $(\G, d)$ is {\it roughly geodesic} if it is $C$-roughly geodesic for some $C \ge 0$.

Let $\Dc_\G$ be the set of metrics $d$ on $\G$ such that $d$ is left-invariant, hyperbolic and quasi-isometric to some (equivalently, any) word metric on $\G$.
Note that every metric $d \in \Dc_\G$ is proper since it is quasi-isometric to a proper metric.
Although $d \in \Dc_\G$ is not necessarily geodesic, it is roughly geodesic by a result of Bonk and Schramm \cite[Proposition 5.6]{BonkSchramm}.
We repeatedly use the following lemma called the Morse lemma when we change metrics in $\Dc_\G$.

\begin{lemma}\label{Lem:Morse}
Let $(\G, d)$ be a proper $C$-roughly geodesic $\d$-hyperbolic space.
For any $(L, K)$-quasi-geodesic $\phi$ in $(\G, d)$, there exists a $C$-rough geodesic $\phi'$ such that $\phi$ and $\phi'$ are within Hausdorff distance $D$ where $D$ depends only on $C, L, K$ and the hyperbolic constant $\d$ of $(\G, d)$.
\end{lemma}

See the proof when $d$ is geodesic \cite[Th\'eor\`emes 21 et 25, Chapitre 5]{GhysdelaHarpe};
it is adapted when $d$ is roughly geodesic, cf.\ \cite[proof of Proposition 5.6]{BonkSchramm}.

\subsubsection{Boundary at infinity}

Let us denote by $\partial (\G, d)$ the {\it geometric boundary} (or {\it boundary}) of $(\G, d)$.
This is a compact metrizable space consisting of equivalence classes of divergent sequences in $(\G, d)$.
Recall that a sequence $\{x_n\}_{n=0}^\infty$ in $(\G, d)$ is {\it divergent} if $(x_n|x_m)_w \to \infty$ as $n, m \to \infty$ for some (equivalently, any) $w \in \G$.
Two divergent sequences $\{x_n\}_{n=0}^\infty$ and $\{y_n\}_{n=0}^\infty$ are equivalent $\{x_n\}_{n=0}^\infty \sim \{y_n\}_{n=0}^\infty$
if $(x_n|y_m)_w \to \infty$ as $n, m \to \infty$.
The boundary $\partial (\G, d)$ is the set of equivalence classes of divergent sequences in $(\G, d)$.
For any $d \in \Dc_\G$, the corresponding boundaries $\partial (\G, d)$ are all homeomorphic each other.
We denote by $\partial \G$ the underlying topological space of the boundary for $d \in \Dc_\G$.
The group $\G$ acts on $\partial \G$ continuously by left multiplications $\G \times \partial \G \to \partial \G$, $(x, \x)\mapsto x\cdot \x$.

Since we work with a roughly geodesic hyperbolic metric $d$ in $\G$,
we record the following lemma which says that there exists a constant $C \ge 0$ such that any two distinct points in the boundary $\partial \G$ are joined by a $C$-rough geodesic in $(\G, d)$.

\begin{lemma}\cite[Proposition 5.2 (3)]{BonkSchramm}\label{Lem:BonkSchramm}
If $(\G, d)$ is a $C$-roughly geodesic $\d$-hyperbolic space,
then there exists a constant $C'=C'(\d, C) \ge 0$ such that
for any two points $\x, \y \in \partial \G$ with $\x \neq \y$, there is a $C'$-rough geodesic $\phi:\R \to \G$ satisfying that
$\phi(-t) \to \x$ and $\phi(t) \to \y$ as $t \to \infty$ in $\G \cup \partial \G$.
\end{lemma}

\subsubsection{Gauges in the boundary}\label{Sec:guage}
Fix a base point $o$ corresponding to the identity element of $\G$.
We extend the Gromov product relative to $d \in \Dc_\G$ to $\G \cup \partial \G$ by setting
\[
(\x|\y)_o:=\sup \Big\{\liminf_{n \to \infty}(x_n|y_n)_o \ : \ \{x_n\}_{n=0}^\infty \in \x, \ \{y_n\}_{n=0}^\infty \in \y\Big\}.
\]
If $d \in \Dc_\G$ is $\d$-hyperbolic, then
for any two pairs of equivalent sequences $\{x_n\}_{n=0}^\infty \sim \{x_n'\}_{n=0}^\infty$ and $\{y_n\}_{n=0}^\infty \sim \{y_n'\}_{n=0}^\infty$, 
\[
\liminf_{n \to \infty}(x_n'|y_n')_o \ge \limsup_{n \to \infty}(x_n|y_n)_o-2\d.
\]
Hence if $\{x_n\}_{n=0}^\infty$ and $\{y_n\}_{n=0}^\infty$ converge to $\x$ and $\y$, respectively in $\G \cup \partial \G$,
then
\begin{equation}\label{Eq:G-prod}
(\x|\y)_o-2\d \le \liminf_{n \to \infty}(x_n|y_n)_o \le \limsup_{n \to \infty}(x_n|y_n)_o \le (\x|\y)_o+2 \d,
\end{equation}
and we have the $\d$-hyperbolic inequality on $\G \cup \partial \G$,
\begin{equation}\label{Eq:delta}
(\x|\y)_o \ge \min\Big\{(\x|\z)_o, (\z|\y)_o\Big\}-3\d \quad \text{for $\x, \y, \z \in \G \cup \partial \G$}.
\end{equation}
For $d \in \Dc_\G$, let us define
\[
\rho(\x, \y):=\exp\(-\(\x|\y\)_o\) \quad \text{for $\x, \y \in \partial \G$}.
\]
The $\rho$ defines a {\it quasi-metric} in $\partial \G$, i.e., $\rho$ satisfies that $\rho(\x, \y)=0$ if and only if $\x=\y$, and $\rho(\x, \y)=\rho(\y, \x)$, 
moreover there exists a constant $C \ge 1$ such that
\[
\rho(\x, \z) \le C\(\rho(\x, \y)+\rho(\y, \z)\) \quad \text{for $\x, \y, \z \in \partial \G$}.
\]
In fact, there exists an $\e \in (0, 1)$ such that $\rho^\e$ is bi-Lipschitz equivalent to a genuine metric on $\partial \G$ (e.g., \cite[Proposition 10, Section 3, Chapitre 7]{GhysdelaHarpe}).
However, in order to avoid introducing an additional parameter $\e$, we mainly use $\rho$ in the boundary $\partial \G$.
Let $\partial^2 \G:=\{(\x, \y) \in (\partial \G)^2 \ : \ \x \neq \y \}$.
We define $\rho_\times$ by
\[
\rho_\times((\x_1, \y_1), (\x_2, \y_2)):=\max\{\rho(\x_1, \x_2), \rho(\y_1, \y_2)\} \quad \text{for $(\x_1, \y_1), (\x_2, \y_2) \in \partial^2 \G$}.
\]
For $d \in \Dc_\G$,
we call $\rho$ and $\rho_\times$ {\it gauges} associated with $d$ on $\partial \G$ and $\partial^2 \G$, respectively.

\subsubsection{Shadows}\label{Sec:shadows}
For $d \in \Dc_\G$, we define shadows on $\partial \G$; they behave similarly to balls relative to $\rho$, but it is better suited to control measures on the boundary.
Moreover, shadows are less sensitive to the change of metrics $d$ in $\Dc_\G$ while $\rho$ vitally depends on $d$.
Let $(\G, d)$ be a $C$-roughly geodesic hyperbolic space and fix a point $o \in \G$.
For any $x \in \G$ and any $R \ge 0$, let us define $O_o(x, R)$ the {\it shadow from $o$} as
the set of points $\x \in \partial \G$ such that some $C$-rough geodesic ray from $o$ converging to $\x$ intersects 
$B_d(x, R):=\{y \in \G \ : \ d(x, y) \le R\}$.
The next lemma follows from the definitions of shadows and the gauge $\rho$ by the $\d$-hyperbolicity.

\begin{lemma}[cf.\ Proposition 2.1 in \cite{BHM11}]\label{Lem:shadows-balls}
Fix $d \in \Dc_\G$.
Then there exist constants $R_0, C>0$ such that for all $R \ge R_0$,
for all $x \in \G$ and all $\x \in \partial \G$, we have that
\[
B_\rho(\x, (1/C)e^{-d(o, x)+R}) \subset O_o(x, R) \subset B_\rho(\x, Ce^{-d(o, x)+R}),
\]
where $B_\rho(\x, r):=\{\y \in \partial \G \ : \ \rho(\x, \y) \le r\}$.
\end{lemma}

For $d, d' \in \Dc_\G$, let $O_o(x, R)$ and $O'_o(x, R)$ be shadows defined in $(\G, d)$ and $(\G, d')$, respectively.
Since $d$ and $d'$ are quasi-isometric, any $C$-rough geodesic $\phi$ in $(\G, d)$ is a $(L, K)$-quasi-geodesic in $(\G, d')$, and thus there exists a $C'$-geodesic in $\phi'$ in $(\G, d')$ such that $\phi$ and $\phi'$ are within Hausdorff distance at most $D$ in $(\G, d')$, where the constant $D$ depends only on the constants involving by Lemma \ref{Lem:Morse}. 
Therefore we have $O_o(x, R) \subset O'_o(x, R')$ for $R'=LR+K+D$.
Note that shadows around the same point $x$ but different metrics in $\Dc_\G$ are comparable up to changing the thickness of shadows independent of $x$.

\subsection{Hausdorff dimension of measures}\label{Sec:dimension}

Let $(X, \rho)$ be a space endowed with a gauge $\rho$.
Examples we have in mind are $(\partial \G, \rho)$ and $(\partial^2 \G, \rho_\times)$.
If $\rho^\e$ is bi-Lipschitz equivalent to a genuine metric $d_\e$ for some $0<\e<1$, then the Hausdorff dimension $\dim_H(A, \rho)$ of a set $A$ relative to the gauge $\rho$ is 
$\e \cdot \dim_H(A, d_\e)$.

For any subset $E$ of $X$, letting
$\rho(E):=\sup\{\rho(x, y) \ : \ x, y \in E \}$,
we define for every $D \ge 0$ and $\D>0$,
\[
\Hc_\D^D(E, \rho):=\inf\Big\{\sum_{i=0}^\infty \rho(E_i)^D \ : \ \text{$E \subset \bigcup_{i=0}^\infty E_i$ and $\rho(E_i) \le \D$}\Big\}.
\]
The $D$-\textit{Hausdorff measure} of a set $E$ is defined by
\[
\Hc^D(E, \rho):=\sup_{\D>0}\Hc_\D^D(E, \rho)=\lim_{\D \to 0}\Hc_\D^D(E, \rho).
\]
The \textit{Hausdorff dimension} of a set $E$ in $(X, \rho)$ is defined by
\[
\dim_H(E, \rho):=\inf\Big\{D \ge 0 \ : \ \Hc^D(E, \rho)=0\Big\}=\sup\Big\{D \ge 0 \ : \ \Hc^D(E, \rho)>0 \Big\}.
\]

\begin{definition}
Let $\n$ be a Borel measure on $(X, \rho)$.
We define the {\it lower Hausdorff dimension} of $\n$ by
\[
\underline{\dim}_H(\n, \rho):=\inf\Big\{\dim_H(E, \rho) \ : \ \n(E)>0,\ \text{$E$ is Borel}\Big\},
\]
and the {\it upper Hausdorff dimension} of $\n$ by
\[
\overline{\dim}_H(\n, \rho):=\inf\Big\{\dim_H(E, \rho) \ : \ \n(X \setminus E)=0,\ \text{$E$ is Borel}\Big\}.
\]
If we have $\underline{\dim}_H(\n, \rho)=\overline{\dim}_H(\n, \rho)$, then we say that the value is the {\it Hausdorff dimension} of $\n$ and denote it by
$\dim_H(\n, \rho)$. 
\end{definition}

In order to estimate the Hausdorff dimension of a measure, we use the following Frostman-type lemma, which relates the dimension of $\n$ to the pointwise behavior $\n\(B_\rho(\x, r)\)$ as $r \to 0$ at each point $\x \in X$.

\begin{lemma}[Cf.\ Sect.\ 8.7 in \cite{Heinonen}]\label{Lem:Frostman}
Let $\n$ be a finite Borel measure on $X$.
If there exist $D_1, D_2 \ge 0$ such that
\[
D_1 \le \liminf_{r \to 0}\frac{\log \n\(B_\rho(\x, r)\)}{\log r} \le D_2 \quad \text{for $\n$-almost every $\x$ in $X$},
\]
where $B_\rho(\x, r)=\{\y \in X \ : \ \rho(\x, \y) \le r\}$,
then
$D_1 \le \underline{\dim}_H(\n, \rho) \le \overline{\dim}_H(\n, \rho) \le D_2$.
\end{lemma}
 
 Note that by Lemma \ref{Lem:Frostman}, for every finite Borel measure $\n$ on $(X, \rho)$, we have
 \[
 \underline{\dim}_H(\n, \rho)=\inf_{\text{$\n$}}\,\liminf_{r \to 0}\frac{\log \n\(B_\rho(\x, r)\)}{\log r}
 \quad \text{and} \quad
 \overline{\dim}_H(\n, \rho)=\sup_{\text{$\n$}}\,\liminf_{r \to 0}\frac{\log \n\(B_\rho(\x, r)\)}{\log r}
 \]
 where the infimum (resp.\ supremum) stands for the essential infimum (resp.\ supremum) relative to $\n$.
 If $\n$ is a possibly infinite but $\s$-finite Borel measure on $(X, \rho)$, then taking compact sets $X_n$ of $X$ for $n=0, 1, \dots$ such that $X=\bigcup_{n=0}^\infty X_n$, and letting $\n|_{X_n}:=\n(\ \cdot \ \cap X_n)$, 
 we have
 \[
 \underline{\dim}_H(\n, \rho)=\inf_{n=0, 1, \dots}\underline{\dim}_H(\n|_{X_n}, \rho) \quad \text{and} \quad \overline{\dim}_H(\n, \rho)=\sup_{n=0, 1, \dots}\overline{\dim}_H(\n|_{X_n}, \rho).
 \]
 We will apply these facts to Radon measures on $(\partial^2 \G, \rho_\times)$ in oder to estimate their Hausdorff dimensions.

\subsection{Busemann (quasi-)cocycles}\label{Sec:Busemann}

For any $d \in \Dc_\G$ and for $w \in \G$, let us define the Busemann function based at $w$ by
\[
b_w:\G \times \partial \G \to \R, \quad b_w(x, \x):=\sup\Big\{\limsup_{n \to \infty}\(d(x, z_n)-d(w, z_n)\) \ : \ \{z_n\}_{n=0}^\infty \in \x \Big\}.
\]
We focus on the Busemann function $b_o$ based at $o$.
Noting the identity
\begin{align*}
d(x, z)-d(o, z)=d(o, x)-2(x|z)_o,
\end{align*}
by \eqref{Eq:G-prod}, 
we have that
$\Big|b_o(x, \x)-\(d(o, x)-2(x|\x)_o\)\Big| \le 4\d$ for $(x, \x) \in \G \times \partial \G$.
This shows that for $x \in \G$ and $\x \in O_o(x, R)$,
\[
-d(o, x) \le b_o(x, \x) \le -d(o, x)+C_{R, \d}, 
\]
where in fact the first inequality holds for all $\x \in \partial \G$ and all $x \in \G$ by the triangle inequality.
Moreover, $b_o$ satisfies the (cocycle) identity up to an additive constant $4\d$,
\begin{equation}\label{Eq:B-cocycle}
\Big|b_o(xy, \x)-\(b_o(y, x^{-1}\x)+b_o(x, \x)\)\Big| \le 4\d \quad \text{for $x, y \in \G, \x \in \partial \G$}.
\end{equation}
Combining the definition of the Gromov product on $\G \cup \partial \G$, we have
\[
b_o(x, \x)+b_o(x, \y)=2(\x|\y)_x-2(\x|\y)_o+O_\d \quad \text{and} \quad b_o(x, \x)-b_o(x, \y)=-2(\x|x)_o+2(\y|x)_o+O_\d,
\]
for $x \in \G$ and two distinct points $\x, \y \in \partial \G$.

\subsection{Strongly hyperbolic metrics}\label{Sec:strongly_hyperbolic}

Let us introduce a special class of metrics which behave regularly at infinity, following \cite{NicaSpakula}.

\begin{definition}\label{Def:strongly_hyperbolic}
We say that a hyperbolic metric $d$ on $\G$ is {\it strongly hyperbolic} if 
there exist constants $L \ge 0$, $s>0$ and $R_0\ge 0$ such that for all $x, x', y, y' \in \G$, and all $R \ge R_0$,
if $d(x, y)-d(x, x')-d(y, y')+d(x', y') \ge R$,
then
\[
|d(x, y)-d(x', y)-d(x, y')+d(x', y')| \le Le^{-s R}.
\] 
\end{definition}

Nica and {\v S}pakula have shown that $d$ is strongly hyperbolic if and only if
there exists a constant $\e>0$ such that
\[
\exp\(-\e(x|y)_o\) \le \exp\(-\e(x|z)_o\)+\exp\(-\e\(z|y\)_o\),
\]
for all $x, y, z, o \in \G$ \cite[Lemma 6.2, Definition 4.1]{NicaSpakula}.
In this form of the definition, it is more transparent to see that the strongly hyperbolicity implies that the (usual) hyperbolicity and the corresponding Gromov product extends to $\G \cup \partial \G$ continuously \cite[Theorem 4.2]{NicaSpakula}.

This property of strong hyperbolicity is actually much stronger than what we expect from word metrics (except for very special cases such as word metrics on finite rank free groups with the standard set of generators).
We are interested in strongly hyperbolic metrics which are quasi-isometric to a word metric.
It is known that every non-elementary hyperbolic group admits a strongly hyperbolic metric which is left-invariant and quasi-isometric to a word metric, i.e., there is a strongly hyperbolic metric in $\Dc_\G$.
In this case, since $d$ is roughly geodesic and quasi-isometric to any word metric, 
the strong hyperbolicity of $d$ is equivalent to the following:
for any Cayley graph $\Cay(\G, S)$
there exist constants $L \ge 0$, $s>0$ and $C \ge 0$ depending only on $\Cay(\G, S)$ such that for any $x, x', y, y' \in \G$ and any $n \ge 0$,
if two geodesic segments connecting $x$ and $y$, and $x'$ and $y'$, respectively, have a common geodesic segment of length $n$ in each $C$-neighborhood, then
\begin{equation}\label{Eq:strongly_hyperbolic_word}
|d(x, y)-d(x', y)-d(x, y')+d(x', y')| \le Le^{-s n},
\end{equation}
where $\d$ is a hyperbolic constant of $\Cay(\G, S)$.

\def\CAT{{\rm CAT}}
\begin{example}\label{Ex:hat}
\item[(1)] The hat metric $\wh d$ on $\G$ introduced by Mineyev is a hyperbolic metric which is strongly hyperbolic, left-invariant and quasi-isometric to a word metric \cite[Theorem 32]{MineyevFlow}.
\item[(2)] $\CAT(-1)$ spaces are strongly hyperbolic \cite[Theorem 5.1]{NicaSpakula}, e.g., Riemannian manifolds with sectional curvature at most $-1$.
If a hyperbolic group $\G$ acts on a $\CAT(-1)$ space isometrically and the action is properly discontinuous, cocompact and {\it free},
then for any base point $o$ in the space, the metric $d_o$ induced from the orbit $d_o(x, y):=d(xo, yo)$ is strongly hyperbolic, left-invariant and quasi-isometric to a word metric.
If the action is not free and the stabilizer at $o$ is nontrivial, then $d_o$ is not a genuine metric on $\G$, but up to changing in the rough isometry class one obtains a metric by defining $d_o'(x, y):=d(xo, yo)+\e$ if $x \neq y$ and $d_o'(x, y):=0$ if $x=y$ for a fixed $\e>0$.
Then the metric $d'_o$ on $\G$ is strongly hyperbolic, left-invariant and quasi-isometric to a word metric.
\item[(3)] The Green metric $d_G$ associated with a $\m$-random walk on $\G$ is a hyperbolic metric which is left-invariant and quasi-isometric to a word metric if $\m$ is a finitely supported and symmetric probability measure on $\G$ such that the support generates the group $\G$ as a semigroup \cite[the first part of Corollary 1.2]{BHM11}.
Moreover, $d_G$ is strongly hyperbolic, which follows from a special case of \cite[Theorem 2.9]{GouezelLocalLimit} and \cite[Section 3]{INO} 
(cf.\ \cite[Theorem 6.1]{NicaSpakula}).
\end{example}

We denote by $\wh d$ a hyperbolic metric which is strongly hyperbolic, left-invariant and quasi-isometric to a word metric on $\G$.
This notation is indebted to the hat metric, but $\wh d$ is not intended to the particular metric in Example \ref{Ex:hat} (1).
The property we use is that $\wh d \in \Dc_\G$ and $\wh d$ is strongly hyperbolic in Definition \ref{Def:strongly_hyperbolic}.

Since $\wh d \in \Dc_\G$, the following is a consequence of the Morse lemma for roughly geodesic hyperbolic spaces (Lemma \ref{Lem:Morse}): there exists a constant $C \ge 0$ such that for any $x, y \in \Cay(\G, S)$ and any $z$ on a geodesic connecting $x$ and $y$ in $\Cay(\G, S)$,
we have
\begin{equation}\label{Eq:Ancona}
|\wh d(x, z)+\wh d(z, y) - \wh d(x, y)| \le C.
\end{equation}
Intending to clarify and simplify notations,
we denote the Gromov product by $(\cdot, \cdot)_o$ relative to a general $d \in \Dc_\G$, and by $\lbr \cdot, \cdot \rbr_o$ relative to $\wh d$, respectively.
An important consequence of the strong hyperbolicity is that the Gromov product $\lbr \cdot, \cdot \rbr_o$ extends continuously on $\G \cup \partial \G$.

For the metric $\wh d$, the corresponding Busemann function $\wh b_o$ (based at $o$) is defined as the limit
\begin{equation}\label{Eq:hat-B}
\wh b_o(x, \x)=\lim_{n \to \infty}\(\wh d(x, z_n) -\wh d(o, z_n)\)=\wh d(o, x)-2\lbr x | \x\rbr_o,
\end{equation}
for any sequence $\{z_n\}_{n=0}^\infty$ converging to $\x$.
Note that for any $x \in \G$, the Busemann function $\wh b_o(x, \cdot)$ is continuous on $\partial \G$.
Furthermore,
$\wh b$ satisfies the genuine cocycle identity in \eqref{Eq:B-cocycle},
\[
\wh b_o(xy, \x)=\wh b_o(y, x^{-1}\x)+\wh b_o(x, \x) \quad \text{for $x, y \in \G, \x \in \partial \G$}.
\]
Note that we have the genuine identity
\begin{align}\label{Eq:B-double}
\wh b_o(x, \x)+\wh b_o(x, \y)=2\lbr \x | \y\rbr_{x}-2\lbr \x|\y\rbr_o,
\end{align}
for $x \in \G$ and for any two distinct $\x, \y \in \partial \G$.

\subsection{Patterson-Sullivan measures}

For a metric $d \in \Dc_\G$, let us denote balls relative to $d$ by
$B_d(x, R):=\{y \in \G \ : d(x, y) \le R\}$ for $x \in \G$ and $R \ge 0$.
We define the exponential volume growth rate of $(\G, d)$ by
\[
\gr(d):=\limsup_{n \to \infty}\frac{1}{n}\log |B_d(o, n)|,
\] 
where $|\cdot|$ stands for the cardinality of the set.
We have that $0<\gr(d)<\infty$ as soon as $\G$ is non-elementary (i.e., non-amenable) and $d \in \Dc_\G$.
The Patterson-Sullivan construction yields a probability measure $\m_d$ associated with $d$ on $\partial \G$ satisfying that
for a constant $C>0$,
\begin{equation}\label{Eq:PS}
\frac{1}{C} \exp\(-\gr(d) \, b_o(x, \x)\) \le \frac{d x_\ast \m_d}{d\m_d}(\x) \le C \exp\(-\gr(d)\,  b_o(x, \x)\),
\end{equation}
for all $x \in \G$ and $\m_d$-almost every $\xi \in \partial \G$. 
For the Patterson-Sullivan construction, see \cite[Th\'eor\`eme 5.4]{Coornaert1993} and for the statement in a more general setting, see \cite[Theorem 2.7]{BHM11} which covers the case for $d \in \Dc_\G$.
We call $\m_d$ a {\it Patterson-Sullivan measure} for $d \in \Dc_\G$.
The following is a consequence of \eqref{Eq:PS} (for which we do not reproduce the proof).

\begin{proposition}\label{Prop:shadow}
Fix a large enough constant $R>0$.
Then there exists a constant $C_R \ge 1$ such that for all $x \in \G$,
\[
\frac{1}{C_R}\exp\(-\gr(d)\, d(o, x)\) \le \m_d\(O_o(x, R)\) \le C_R \exp\(-\gr(d)\, d(o, x)\).
\]
\end{proposition}

We use the following doubling property of Patterson-Sullivan measures later.

\begin{lemma}\label{Lem:doubling}
Let us fix a metric $d \in \Dc_\G$.
Then the corresponding Patterson-Sullivan measure $\m_{d}$ is doubling with respect to $\rho$, i.e.,
there exists a constant $C \ge 1$ such that
for any $\x \in \partial \G$ and for any $r>0$,
\[
\m_{d}\(B_\rho(\x, 2r)\) \le C\cdot \m_{d}\(B_\rho(\x, r)\),
\]
where $B_\rho(\x, r)=\{\y \in \partial \G \ : \ \rho(\x, \y) < r\}$.
\end{lemma}
\proof
Proposition \ref{Prop:shadow} shows that a large enough $R>0$, the measures $\m_d(O_o(x, 2R))$ and $\m_d(O_o(x, R))$ are comparable with uniform constants independent of $x \in \G$.
Applying this estimate repeatedly, Lemma \ref{Lem:shadows-balls} yields the claim.
\qed

It is known that the size of (thickened) spheres is comparable with an exponential function.
The following is a version in the form of a weighted sum relative to a word metric.

\begin{lemma}\label{Lem:regular-growth}
For any hyperbolic metric $d$ in $\Dc_\G$ and any word metric $d_S$ on $\G$,
there exist constants $c_1, c_2>0$ such that
\[
c_1\le \sum_{x \in \bS_n}\exp\(-\gr(d)\, d(o, x)\)\le c_2 \quad \text{for all $n \ge 0$},
\]
where $\bS_n:=\{x \in \G \ : \ d_S(o, x)=n \}$ for integers $n \ge 0$.
\end{lemma}

\proof
Let us denote by $O_o(x, R)$ shadows relative to $d$.
Fix a large enough $R>0$.
For any integer $n$, shadows $O_o(x, R)$ for $|x|_S=n$ cover the boundary $\partial \G$, and thus Proposition \ref{Prop:shadow} implies that
\[
1=\m_d(\partial \G) \le \sum_{x \in \bS_n}\m_d(O_o(x, R)) \le C_R \sum_{x \in \bS_n}\exp(-\gr(d)d(o, x)).
\]
We denote by $O_o^w(x, R)$ shadows relative to a word metric.
For any $R>0$, there exists a constant $R'>0$ such that
$O_o(x, R) \subset O_o^w(x, R')$ 
for all $x \in \G$ (cf.\ Section \ref{Sec:shadows}).
Note that $O_o^w(x, R')$ for $|x|_S=n$ cover the boundary $\partial \G$ with bounded overlaps, i.e., for an integer $D \ge 1$
each point in $\partial \G$ is included at most $D$ shadows $O_o^w(x, R')$ with $|x|_S=n$.
Therefore by Proposition \ref{Prop:shadow} we have that
\[
(1/C_R)\sum_{x \in \bS_n}\exp\(-\gr(d)d(o, x)\) \le \sum_{x \in \bS_n}\m_d(O_o^w(x, R')) \le D \m_d(\partial \G)=D.
\]
We conclude the claim.
\qed

\begin{proposition}\label{Prop:invariant_measure}
Let $\partial^2 \G:=(\partial \G)^2 \setminus \{(\x, \x) \ : \ \x \in \partial \G\}$.
For every $d \in \Dc_\G$, there exists a $\G$-invariant Radon measure $\L_d$ on $\partial^2 \G$ such that
\[
\L_d =\f(\x, \y) \exp(2\gr(d)\, (\x|\y)_o)\m_d \otimes \m_d,
\]
where there exist constants $c_1, c_2>0$ such that $c_1\le \f(\x, \y) \le c_2$ for all $(\x, \y) \in \partial^2 \G$.
Moreover, $\L_d(U)>0$ for each open set $U$ in $\partial^2 \G$ and $\L_d(K)<\infty$ for each compact set $K$ in $\partial^2 \G$.
\end{proposition}

\proof
Let $\m_d$ be a Patterson-Sullivan measure for $d \in \Dc_\G$.
Restricting the measure $\m_d \otimes \m_d$ on $\partial^2 \G$,
we define the measure on $\partial^2 \G$ by
\[
\n:=\exp(2\gr(d) (\x|\y)_o)\m_d \otimes \m_d.
\]
Since we have for $x \in \G$ and $\x, \y \in \partial \G$,
$b_o(x, \x)+b_o(x, \y)=2(\x|\y)_x-2(\x|\y)_o\pm C_\d$
(Section \ref{Sec:Busemann})  
and by \eqref{Eq:PS}, 
there exist constants $c_1, c_2>0$ such that
\begin{equation*}
c_1 \le \frac{d x_\ast \n}{d \n}(\x, \y) \le c_2 \quad \text{for all $x \in \G$ and $\n$-almost every $(\x, \y) \in \partial^2 \G$}.
\end{equation*}
For $x, y \in \G$, note that
\[
\frac{d(x y)_\ast \n}{d\n}(\x, \y)=\frac{d x_\ast \n}{d \n}(\x, \y)\frac{d y_\ast \n}{d \n}(x^{-1}\x, x^{-1}\y).
\]
Letting $\f(\x, \y):=\sup_{x \in \G}\frac{d x_\ast \n}{d\n}(\x, \y)$, we have $c_1 \le \f(\x, \y) \le c_2$ and
\[
\f(\x, \y)=\frac{d x_\ast\n}{d\n}(\x, \y)\cdot \f(x^{-1}\x, x^{-1}\y).
\]
Hence if we define $\L_d:=\f(\x, \y)\n$, then $\L_d$ is $\G$-invariant.
Moreover, note that each compact set $K$ in $\partial^2 \G$ is included in
$\{(\x, \y) \in \partial^2 \G \ : \ (\x|\y)_o \le C\}$
for some $C >0$, and thus $\n_d(K)<\infty$.
For each open set $U$ in $\partial^2 \G$, we have $\L_d(U)>0$ since $\m_d(B)>0$ for each open set $B$ in $\partial \G$.
\qed

\begin{remark}
For a strongly hyperbolic metric $\wh d$, the corresponding $\G$-invariant measure is given by 
$\L_{\wh d}=\exp(2\gr(\wh d\,)\, \lbr \x|\y\rbr_o)\m_{\wh d} \otimes \m_{\wh d}$, that is, one may take $\f\equiv 1$ in Proposition \ref{Prop:invariant_measure}.
This is because one is able to take the constant $C=1$ in \eqref{Eq:PS} in the Pattersson-Sullivan construction and the Busemann function $\wh b_o(x, \cdot\,)$ defines a genuine cocycle (without an additive constant) (see also \cite[Section 7.1]{NicaLp}).
\end{remark}

\section{Topological flow spaces}\label{Sec:topological_flow_space}

\subsection{Construction of a flow space}\label{Sec:flow}

Let $\partial^2 \G:=(\partial \G)^2\setminus \{(\x, \x) \ : \ \x \in \partial \G\}$, i.e., the space of ordered pairs of distinct points in $\partial \G$.
We consider the diagonal action of $\G$ on $\partial^2 \G$.
The space $\partial^2 \G \times \R$ is regarded as the set of directed lines with distinct extreme points in $\partial \G$.
Here we give a geometric realization of $\partial^2 \G \times \R$ by rough geodesics on $(\G, \wh d \,)$.
For every $(\x, \y) \in \partial^2 \G$, there exists a $C$-rough geodesic $\phi_{\x, \y}:\R \to (\G, \wh d\,)$ from $\x$ to $\y$ by Lemma \ref{Lem:BonkSchramm}.
We normalize the parametrization of $\phi_{\x, \y}(t)$ by shifting $t \mapsto t+T$ for some $T$ if necessary such that
\[
\wh d(\phi_{\x, \y}(0), o)=\min_{t \in \R} \wh d(\phi_{\x, \y}(t), o).
\]
Then we define the map 
\[
\partial^2 \G \times \R \to (\G, \wh d\,), \qquad (\x, \y, t) \mapsto \phi_{\x, \y}(t).
\]
This map is not canonical; it depends on choices for rough geodesics.
However, note that for any choices $\phi_{\x, \y}$ and $\phi'_{\x, \y}$ for $(\x, \y) \in \partial^2 \G$,
we have $\sup_{t \in \R}\wh d(\phi_{\x, \y}(t), \phi'_{\x, \y}(t)) \le C$
for a constant $C$ depending only on $(\G, \wh d\,)$ by the stability of $C$-rough geodesics.

\begin{lemma}\label{Lem:contraction}
For the metric $\wh d \in \Dc_\G$, let $\wh \rho$ be the corresponding gauge in $\partial \G$.
For any compact set $K$ in $\partial^2 \G$, there exists a constant $c_K>0$ such that
for any $\x, \y, \y' \in \partial \G$ and $x \in \G$, 
if 
$(\x, \y), (\x, \y') \in K$ and $T=-\wh b_o(x^{-1}, \x)$, 
then
\[
\wh \rho(x\y, x\y') \le e^{2c_K-T}\cdot \wh \rho(\y, \y').
\]
\end{lemma}
\proof
Since the Gromov product $\lbr \cdot, \cdot \rbr_o$ is continuous on $\partial^2\G$ and $K$ is a compact subset in $\partial^2 \G$,
there exists a positive constant $c_K$ depending only on $K$ such that
$\lbr \x | \y\rbr_o, \lbr \x | \y'\rbr_o \le c_K$.
Hence 
\[
\wh b_o(x^{-1}, \x)+\wh b_o(x^{-1},\y) \ge -2 c_K \quad \text{and} \quad \wh b_o(x^{-1}, \x)+\wh b_o(x^{-1}, \y') \ge -2 c_K,
\]
and $\wh b_o(x, \y)+\wh b_o(x, \y')=2\lbr \y |\y'\rbr_{x}-2\lbr \y|\y'\rbr_o$ by \eqref{Eq:B-double},
we obtain
\[
\lbr \y|\y'\rbr_{x^{-1}}-\lbr \y|\y'\rbr_o \ge -\wh b_o(x^{-1}, \x)-2c_K.
\]
Therefore, letting $T=-\wh b_o(x^{-1}, \x)$, we have
\[
\wh \rho(x\y, x\y') =e^{-\lbr \y|\y'\rbr_{x^{-1}}}\le e^{2 c_K}e^{-T-\lbr \y|\y'\rbr_o}=e^{2 c_K}e^{-T}\wh \rho(\y, \y'),
\]
as desired.
\qed

Let us define the cocycle $\a:\G \times \partial^2 \G \to \R$ by
\[
\a(x, \x, \y):=\frac{1}{2}\Big(\wh b_o(x^{-1}, \x)-\wh b_o(x^{-1}, \y)\Big) \quad \text{for $x \in \G$, $(\x, \y) \in \partial^2\G$}.
\]
Note that we have
\begin{equation}\label{Eq:alpha}
\a(x, \x, \y)=\lbr x^{-1}|\y \rbr_o-\lbr x^{-1}|\x \rbr_o.
\end{equation}
Then, we define the action of $\G$ on $\partial^2 \G \times \R$ by for $x \in \G$ and for $(\x, \y, t) \in \partial^2 \G \times \R$,
\begin{equation*}\label{Eq:alpha-action}
x\cdot \Big(\x, \y, t\Big):=\Big(x\x, x\y, t-\a(x, \x, \y)\Big).
\end{equation*}
The cocycle identity $\a\(x y, \x, \y\)=\a\(x, y\x, y\y\)+\a\(y, \x, \y\)$ shows that it defines indeed a $\G$-action on $\partial^2 \G \times \R$.
We call this $\G$-action on $\partial^2 \G \times \R$ via $\a$  
the {\it $(\G, \a)$-action}.

\begin{lemma}\label{Lem:Gamma-action}
The $(\G, \a)$-action of $\G$ on $\partial^2 \G \times \R$ is properly discontinuous and cocompact, i.e., $\G \backslash (\partial^2 \G \times \R)$ is compact.
\end{lemma}

\proof
For $N, T>0$, let
\[
F_{N, T}:=\Big\{(\x, \y, t) \in \partial^2 \G \times \R \ : \ |t| \le T, \ \wh \rho(\x, \y) \ge \exp(-N) \Big\}.
\]
First we shall show that for each $N, T>0$, the number of $x$ in $\G$ such that $x F_{N, T} \cap F_{N, T}\neq \emptyset$ is finite.
Suppose that for $x \in \G$ there exists $(\x, \y, t) \in F_{N, T}$ such that $x(\x, \y, t) \in F_{N, T}$.
Then, we have
\[
\lbr \x|\y\rbr_o \le N \quad \text{and} \quad \lbr \x|\y\rbr_{x^{-1}}=\lbr x \x|x\y\rbr_o \le N.
\]
Furthermore, since $|t| \le T$ and $|t-\a(x, \x, \y)| \le T$,
we have by \eqref{Eq:alpha},
\[
|\lbr x^{-1}|\y \rbr_o-\lbr x^{-1}|\x \rbr_o|=|\a(x, \x, \y)| \le 2T,
\]
and
\begin{align*}
\lbr x^{-1}|\x \rbr_o+\lbr x^{-1}|\y\rbr_o	&\le 2\min\{\lbr x^{-1}|\x \rbr_o, \lbr x^{-1}|\y\rbr_o\}+|\lbr x^{-1}|\y \rbr_o-\lbr x^{-1}|\x \rbr_o|\\
									&\le 2 \lbr \x|\y\rbr_o+6\d+2T,
\end{align*}
where the last inequality follows from the $\d$-hyperbolic inequality \eqref{Eq:delta}.
Therefore by \eqref{Eq:hat-B} 
we have
\begin{align*}
\wh d(o, x)	&=\frac{1}{2}\(\wh b_o(x^{-1}, \x)+\wh b_o(x^{-1}, \y)\)+\lbr x^{-1}|\x\rbr_o+\lbr x^{-1}|\y\rbr_o\\
			&=\lbr \x|\y\rbr_{x^{-1}}-\lbr \x|\y\rbr_o+\lbr x^{-1}|\x\rbr_o+\lbr x^{-1}|\y\rbr_o\\
			&\le \lbr \x|\y\rbr_{x^{-1}}+\lbr \x|\y\rbr_o+6\d+2T \le 2N+6\d+2T.
\end{align*}
This shows that $x$ is included in a ball of a finite radius in the proper metric space $(\G, \wh d\,)$, and thus the number of such $x$ is finite.
Since any compact set $K$ in $\partial^2 \G \times \R$ is contained in $F_{N, T}$ for some $N$ and $T$, 
we conclude that the action of $\G$ on $\partial^2 \G \times \R$ is properly discontinuous.

Second we shall show that for some $N$ and $T$, every $\G$-orbit intersects $F_{N, T}$.
For any given $(\x, \y, t) \in \partial^2 \G \times \R$,
take a $C_\d$-rough geodesic $\phi_{\x, \y}$ from $\x$ to $\y$ such that $\phi_{\x, \y}(0)$ is one of the nearest point on $\phi_{\x, \y}$ to $o$.
Then we have 
$\wh b_o(\phi_{\x, \y}(0), \x)=-\lbr \x|\y \rbr_o+O_\d$ and $\wh b_o(\phi_{\x, \y}(0), \y)=-\lbr \x|\y \rbr_o+O_\d$,
and thus
$\a(\phi_{\x, \y}(0)^{-1}, \x, \y)=O_\d$.
Since for any $r \in \R$,
\[
\wh b_o(\phi_{\x, \y}(r), \x)=r+\wh b_o(\phi_{\x, \y}(0), \x)+O_\d, \quad \wh b_o(\phi_{\x, \y}(r), \y)=-r+\wh b_o(\phi_{\x, \y}(0), \y)+O_\d,
\]
we have that
$\a(\phi_{\x, \y}(r)^{-1}, \x, \y)=r+\a(\phi_{\x, \y}(0)^{-1}, \x, \y)+O_\d=r+O_\d$.
Letting $x:=\phi_{\x, \y}(t)^{-1}$, we obtain
\[
\a(x, \x, \y)=\a(\phi_{\x, \y}(t)^{-1}, \x, \y)=t\pm C'_\d.
\]
Since now the base point $o$ is on the $C_\d$-rough geodesic $x\phi_{\x, \y}$ from $x \x$ to $x \y$,
we have $\lbr x\x|x\y \rbr_o \le 3C_\d+2\d$.
Therefore if $N:=3C_\d+2\d$ and $T:=C'_\d$, then
\[
x\cdot \Big(\x, \y, t\Big) = \Big(x\x, x\y, t-\a(x, \x, \y)\Big) \in F_{N, T}.
\]
The set $F_{N, T}$ is compact, and thus the action $\G$ on $\partial^2 \G\times \R$ is cocompact.
\qed

Let $\Fc_\a:=\G \backslash \(\partial^2 \G \times \R\)$ be the quotient topological space of $\partial^2 \G \times \R$ by the $(\G, \a)$-action.
The space $\Fc_\a$ is compact by Lemma \ref{Lem:Gamma-action}.
We shall define a flow on $\Fc_\a$, namely, a continuous $\R$-action on $\Fc_\a$.

We define the action of $\R$ on $\partial^2 \G \times \R$ by for $t \in \R$, $(\x, \y, s) \in \partial^2 \G \times \R$,
\[
\wt \F_t(\x, \y, s):=(\x, \y, t+s).
\]
Since this $\R$-action and the $(\G, \a)$-action commute on $\partial^2 \G \times \R$,
the action $\{\wt \F_t\}_{t \in \R}$ descends to $\Fc_\a$.
We define the induced flow by for $t \in \R$, $[\x, \y, s] \in \Fc_\a$,
\[
\F_t: \Fc_\a \to \Fc_\a, \quad [\x, \y, s] \mapsto [\x, \y, t+s].
\]
The action of $\R$ via $\F_t$ is continuous on $\Fc_\a$.
We call $\{\F_t\}_{t \in \R}$ a {\it topological flow} on $\Fc_\a$.

\begin{remark}
The space $\Fc_\a$ plays a role of the (total space of) unit tangent bundle over a compact Riemannian manifold where the geodesic flow is defined in the classical sense.
Let us mention a metric structure on $\Fc_\a$ not just the topology although we do not use it for our purpose.
Associated with a metric $\wh d$ in $\G$, 
there is a metric $d_\ast$ in $\partial^2 \G \times \R$ such that the $(\G, \a)$-action on $(\partial^2 \G \times \R, d_\ast)$ is isometric
\cite[Theorem 60 (d)]{MineyevFlow}.
One may define the metric $d_{\Fc_\a}$ in $\Fc_\a$ by
\[
d_{\Fc_\a}([x], [x']):=\inf_{\g \in \G}d_\ast(x, \g x') \quad \text{for $x, x' \in \partial^2 \G \times \R$},
\]
where the infimum is attained since the $(\G, \a)$-action on $\partial^2 \G \times \R$ is cocompact.
\end{remark}

\subsection{Measures on $\partial^2 \G \times \R$ and $\Fc_\a$}

Let $\L$ be a $\G$-invariant Radon measure on $\partial^2\G$.
(Recall that a measure on a locally compact Hausdorff space is Radon if it is finite on every compact set and Borel regular.)
Let $dt$ be the (normalized) Lebesgue measure on $\R$.
Then, the measure $\L \otimes dt$ is Radon on $\partial^2 \G \times \R$.
Furthermore, the measure $\L \otimes dt$ is $(\G, \a)$-invariant as well as $\{\wt \F_t\}_{t \in \R}$-invariant.
For each continuous function with compact support $f \in C_c(\partial^2 \G \times \R)$,
let
\[
\wbar f(\x, \y, t):=\sum_{x \in \G}f\(x\cdot (\x, \y, t)\),
\]
where the summation runs over a finite number of $x$ for each point by Lemma \ref{Lem:Gamma-action}.
Since $\wbar f$ is $\G$-invariant on $\partial^2 \G \times \R$, the function $\wbar f$ can be defined on $\Fc_\a$.
We use the same symbol for the function $\wbar f$ defined on $\Fc_\a$.

\begin{lemma}\label{Lem:quotient-measure}
For any $\G$-invariant Radon measure $\L$ on $\partial^2 \G$, there exists a unique finite Radon measure $m$ on $\Fc_\a$
such that
\[
\int_{\partial^2 \G \times \R}f\,d\L\otimes dt=\int_{\Fc_\a}\wbar f\,dm,
\]
for any compactly supported continuous function $f$ on $\partial^2 \G \times \R$.
Moreover, the measure $m$ is invariant under the topological flow $\{\F_t\}_{t \in \R}$.
\end{lemma}

\proof
Note that for any continuous function $\f$ on $\Fc_\a$, there exists a $f \in C_c(\partial^2 \G \times \R)$ such that
$\wbar f=\f$.
Indeed, let us denote the quotient map by $\pi:\partial^2 \G\times \R \to \Fc_\a$.
There exists a compact set $F$ on $\partial^2 \G \times \R$ such that $\pi(F)=\Fc_\a$ by Lemma \ref{Lem:Gamma-action}.
By the Urysohn lemma, there exists a continuous function $\chi$ with values in $[0, 1]$ on $\partial^2 \G \times \R$ such that $\chi \equiv 1$ on $F$ and $\chi \equiv 0$ outside a relatively compact neighborhood of $F$.
Taking such a function $\chi$,
define
\[
f(x):=\frac{\chi(x)}{\sum_{\g \in \G}\chi(\g x)}\f(\pi(x)).
\]
Then $f$ has a compact support and $\wbar f=\f$. 
Moreover, if $\f \ge 0$, then $f \ge 0$.

For any continuous function $\f$ on $\Fc_\a$, the map
\[
\f \mapsto \int_{\partial^2 \G \times \R}f(x)d\L\otimes dt,
\]
is well-defined; independent of the choice of $f$ such that $\wbar f=\f$, and defines a positive bounded linear functional on the space of continuous functions $C(\Fc_\a)$.
Therefore the Riesz representation theorem yields a unique finite Radon measure $m$.

Since $\L \otimes dt$ is $\{\wt \F_t\}_{t \in \R}$-invariant and the $(\G, \a)$-action and $\{\wt \F_t\}_{t \in \R}$ commute on $\partial^2 \G \times \R$,
for any $\f \in C(\Fc_\a)$ and any $t \in \R$,
\[
\int_{\Fc_\a}\f \circ \F_t\,dm=\int_{\Fc_\a}\f\,dm.
\]
This shows that $m$ is invariant under $\{\F_t\}_{t \in \R}$.
\qed

\begin{remark}\label{Rem:quotient-measure}
If we take a Borel fundamental domain $D$ and a measurable section $\i:\Fc_\a \to D \subset \partial^2 \G \times \R$,
then we have $\L \otimes dt=\sum_{\g \in \G}\g. \i_\ast m$.
If the group $\G$ acts on $\partial^2 \G \times \R$ {\it freely}, then the measure $m$ in Lemma \ref{Lem:quotient-measure} is obtained by the restriction of the measure $\L \otimes dt$ on $D$, i.e., $\L \otimes dt|_{D}=\i_\ast m$.
However, if the $\G$-action is not free, then $\i_\ast m$ is not the restriction of $\L \otimes dt$ on $D$.
See \cite[Section 2.6]{PPS}.
\end{remark}

\subsection{Flow invariant measures}

For any $d \in \Dc_\G$, let $\m_d$ be a Patterson-Sullivan measure on $\partial^2 \G$ and
$\L_d$ be the associated $\G$-invariant Radon measure on $\partial^2 \G$ constructed in Proposition \ref{Prop:invariant_measure}.
Let us consider the $(\G, \a)$-invariant measure $\L_d \otimes dt$ on $\partial^2 \G \times \R$.
Lemma \ref{Lem:quotient-measure} implies that there exists a unique finite Radon measure $m_d$ on $\Fc_\a$ such that
$m_d$ is $\{\F_t\}_{t \in \R}$-invariant and
\[
\int_{\partial^2 \G \times \R}f\,d\L_d\otimes dt=\int_{\Fc_\a}\wbar f\,dm_d \quad \text{for $f \in C_c(\partial^2 \G \times \R)$}.
\]
We normalize $m_d$ with $m_d(\Fc_\a)=1$.

\begin{theorem}\label{Thm:ergodic}
Fix $d \in \Dc_\G$.
Then for any $\f$ in $L^1(\Fc_\a, m_d)$, we have
\[
\lim_{T \to \infty}\frac{1}{T}\int_0^T \f\circ \F_t(x)\,dt=\int_{\Fc_\a}\f\,dm_d, \quad \text{for $m_d$-almost every $x$}.
\]
In particular, for any $d \in \Dc_\G$, the measure $m_d$ is ergodic under the topological flow $\{\F_t\}_{t \in \R}$ on $\Fc_\a$, i.e., for any Borel set $A$ in $\Fc_\a$ with $\F_{-t}(A)=A$ for all $t \in \R$, one has $m_d(A)=0$ or $m_d(\Fc_\a \setminus A)=0$.
\end{theorem}

\proof
For $d \in \Dc_\G$, let us write $\L=\L_d$ and $m=m_d$.
For any $\f \in L^1(\Fc_\a, m)$,
the Birkhoff ergodic theorem implies that for $m$-almost every $x \in \Fc_\a$,
there exists a limit
\[
\f_\infty(x):=\lim_{T \to \infty}\frac{1}{T}\int_0^T \f\circ \F_t(x)\,dt,
\]
and $\f_\infty$ is defined on $m$-almost everywhere on $\Fc_\a$ such that $\f_\infty \circ \F_t=\f_\infty$ for any $t \in \R$.
(This follows by considering a Birkhoff sum of $\int_0^1 \f\circ \F_t(x)\,dt$ and integers $T$.)
The convergence is also in $L^1(\Fc_\a, m)$.
We shall show that $\f_\infty$ is constant $m$-almost everywhere.
Note that if $\f_\infty$ coincides with a constant $m$-almost everywhere, then it has to be $\int_{\Fc_\a}\f\,dm$.
By the Fatou lemma, we have $\|\f_\infty\|_1 \le \|\f\|_1$.
Hence it suffices to show the claim for continuous functions $\f \in C(\Fc_\a)$.

Let $\pi:\partial^2 \G \times \R \to \Fc_\a$ be the quotient map.
Taking a Borel fundamental domain $D$ in $\partial^2 \G\times \R$ and a measurable section $\i:\Fc_\a \to D$, 
we have $\L\otimes dt=\sum_{\g \in \G}\g.\i_\ast m$ (Remark \ref{Rem:quotient-measure}).
Hence for $\L\otimes dt$-almost every $(\x, \y, s) \in \partial^2 \G \times \R$,
\begin{equation}\label{Eq:thm-ergodic}
\lim_{T \to \infty}\frac{1}{T}\int_0^T (\f \circ \pi)\(\wt \F_t(\x, \y, s)\)\,dt=\f_\infty \circ \pi (\x, \y, s).
\end{equation}
In fact, since this holds for $\L\otimes dt$-almost every $(\x, \y, s)$ and $\f$ is continuous on a compact space $\Fc_\a$, the limit in \eqref{Eq:thm-ergodic} exists for $\L$-almost every $(\x, \y) \in \partial^2 \G$ and for all $s \in \R$. 
Let $\wt \f_\infty:=\f_\infty \circ \pi$ on $\partial^2 \G \times \R$.
Since $\f_\infty$ is $\{\F_t\}_{t \in \R}$-invariant, the lift $\wt \f_\infty$ is $\{\wt \F_t\}_{t \in \R}$-invariant.
Therefore $\wt \f_\infty$ is defined on $\partial^2 \G$, and we write $\wt \f_\infty(\x, \y)$ for $(\x, \y) \in \partial^2 \G$.
Noting that $\L$ and $\m \otimes \m$ are mutually absolutely continuous on $\partial^2 \G$,
we have for $\m \otimes \m$-almost every $(\x, \y) \in \partial^2 \G$, and all $s \in \R$,
\begin{equation}\label{Eq:thm-ergodic2}
\lim_{T \to \infty}\frac{1}{T}\int_0^T (\f \circ \pi)\(\wt \F_t(\x, \y, s)\)\,dt=\wt \f_\infty(\x, \y).
\end{equation}
Let $K$ be any compact set in $\partial^2 \G$.
Then there exists a constant $c_K>0$ such that
for any $(\x, \y), (\x, \y') \in K$, $\lbr \x|\y\rbr_o, \lbr \x|\y'\rbr_o \le c_K$.
For $\x \in \partial \G$, if we take a $C_\d$-rough geodesic $\phi$ from $o$ to $\x$, then for any $r\ge 0$,
\[
\wh b_o(\phi(r), \x)=-r\pm2C_\d.
\]
For any $T>0$, let $x:=\phi(T)^{-1}$, then $\wh b_o(x^{-1}, \x)=-T \pm 2C_\d$.
Hence Lemma \ref{Lem:contraction} implies that
\[
\wh \rho(x\y, x\y') \le e^{2c_K+2C_\d}\cdot e^{-T},
\]
and thus for any $\e>0$, there exists $x \in \G$ depending only on $\x$ such that for any $(\x, \y), (\x, \y') \in K$,
we have $\wh \rho(x\y, x\y')\le \e$.
Moreover, since 
\[
\a(x, \x, \y)=\frac{1}{2}\(\wh b_o(x^{-1}, \x)-\wh b_o(x^{-1}, \y)\)=\wh b_o(x^{-1}, \x)-\frac{1}{2}\(\wh b_o(x^{-1}, \x)+\wh b_o(x^{-1}, \y)\),
\]
letting
\[
s_1:=\frac{1}{2}\(\wh b_o(x^{-1}, \x)+\wh b_o(x^{-1}, \y)\) \quad \text{and} \quad s_2:=\frac{1}{2}\(\wh b_o(x^{-1}, \x)+\wh b_o(x^{-1}, \y')\),
\]
we have
$s_1+\a(x, \x, \y)=s_2+\a(x, \x, \y')$.
Since $\f \circ \pi$ is uniformly continuous on $\partial^2 \G \times \R$, 
for any $\e>0$ and for any $(\x, \y), (\x, \y') \in K$,
there exists a $x \in \G$ such that for any $t \in \R$,
\[
|\f \circ \pi(x\cdot (\x, \y, s_1+t))-\f \circ \pi(x\cdot (\x, \y', s_2+t))| \le \e.
\]
Therefore since $\f \circ \pi=\f\circ \pi (x \ \cdot \ )$ for all $x \in \G$, by \eqref{Eq:thm-ergodic2},
we obtain for $\m \otimes \m$-almost every $(\x, \y), (\x, \y') \in K$,
\[
|\wt \f_\infty(\x, \y)-\wt \f_\infty(\x, \y')| \le \e.
\] 
This shows that for $\m\otimes \m$-almost every $(\x, \y), (\x, \y') \in K$, we have $\wt \f_\infty(\x, \y)=\wt \f_\infty(\x, \y')$.
Since $\partial^2 \G$ is $\s$-compact, for $\m\otimes \m$-almost every $(\x, \y), (\x, \y') \in \partial^2 \G$,
we have $\wt \f_\infty(\x, \y)=\wt \f_\infty(\x, \y')$, i.e.,
$\wt \f_\infty(\x, \ \cdot \ )$ is $\m$-almost everywhere constant for $\m$-almost every $\x \in \partial \G$.
Changing the role of $\x$ and $\y$ in $\partial^2 \G$, 
we also show that $\wt \f_\infty(\ \cdot \ , \y)$ is $\m$-almost everywhere constant for $\m$-almost every $\y \in \partial \G$.
We conclude that $\wt \f_\infty$ is $\m\otimes \m$-almost everywhere constant on $\partial^2 \G$, and thus $\L$-almost everywhere constant on $\partial^2 \G$.
Therefore since $\wt \f_\infty=\f_\infty \circ \pi$, the limiting function $\f_\infty$ is $m$-almost everywhere constant and we obtain the claim.
\qed

\begin{corollary}\label{Cor:double-ergodicity}
For any $d \in \Dc_\G$, the corresponding $\G$-invariant Radon measure $\L_d$ on $\partial^2 \G$ is ergodic with respect to the $\G$-action on $\partial^2 \G$.
\end{corollary}
\proof
Let $A$ be any $\G$-invariant Borel set in $\partial^2 \G$ and consider $A \times \R$.
Since $A$ is $\G$-invariant, $A \times \R$ is $(\G, \a)$-invariant on $\partial^2 \G \times \R$ and thus $A \times \R$ defines the set $[A \times \R]$ in $\Fc_\a$. 
The set $[A \times \R]$ is $\{\F_t\}_{t \in \R}$-invariant, and it has $m_d$-measure $0$ or $1$ since $m_d$ is ergodic by Theorem \ref{Thm:ergodic}.
This shows that $A \times \R$ has either null or co-null $\L_d \otimes dt$-measure and $A$ has either null or co-null $\L_d$-measure.
Therefore $\L_d$ is ergodic. 
\qed

\section{Coding via automatic structures}\label{Sec:coding}

\subsection{Automatic structures}

Let $S$ be a finite set of generators such that $S=S^{-1}$ in $\G$.
An {\it automaton} $(\Ac, w, S)$ is a directed graph $\Ac=(V, E, s_\ast)$ where $s_\ast$ is a vertex called the {\it initial state}, together with a {\it labeling} $w:E \to S$ on edges by $S$.
For a directed path $\o=(e_0, e_1, \dots, e_{n-1})$ in the graph $\Ac$ where the terminus of $e_i$ is the origin of $e_{i+1}$,
we associate a path $w(\o)$ in the Cayley graph of $(\G, S)$ issuing from the identity
$o, w(e_0), w(e_0)w(e_1), \dots, w(e_0)\cdots w(e_{n-1})$.
Let $w_\ast(\o)$ be the terminus of the path $w(\o)$, namely, $w_\ast(\o)=w(e_0)\cdots w(e_{n-1})$.

\begin{definition}
We say that an automaton $(\Ac, w, S)$ where $\Ac=(V, E, s_\ast)$ and a labeling $w:E \to S$ is a {\it strongly Markov automatic structure}
if 
\begin{itemize}
\item[(1)] each vertex in $V$ can be reached by a directed path from the initial state $s_\ast$,
\item[(2)] for every directed path $\o$ in $\Ac$, the path $w(\o)$ is a geodesic in the Cayley graph of $(\G, S)$, and
\item[(3)] the map $w_\ast$ assigning the terminus of a path defines a bijection from the set of directed paths from $s_\ast$ in $\Ac$ to $\G$. 
\end{itemize}
\end{definition}

Every hyperbolic group admits a strongly Markov automatic structure for any finite set of generators $S$ with $S=S^{-1}$ \cite{Cannon} and \cite[Section 3.2]{Calegari}.
We fix an automaton $(\Ac, w, S)$ for $(\G, S)$.
Let $\SS^\ast$ be the set of finite directed paths (not necessarily from $s_\ast$) in the graph $\Ac$,
and $\SS^+$ be the set of semi-infinite paths $(e_i)_{n=0}^\infty$.
Let $\wbar \SS^+:=\SS^\ast \cup \SS^+$ be the set of unilateral paths.
We define the map $w_\ast:\wbar \SS^+ \to \G \cup \partial \G$, where
to every path $\o$ in $\wbar \SS^+$, we associate a point in $\G\cup \partial \G$ as the terminus of geodesic segment or a geodesic ray $w(\o)$ starting from the identity $o$ in $\Cay(\G, S)$.

Actually we mainly use the space of bilateral paths based on $\Ac$.
We will work with the space $\wbar \SS^+$ to construct an appropriate measures in Section \ref{Subsec:thermodynamic}.

\subsection{Bilateral paths}
\def\inv{{\rm inv}}
Let us define $\SS$ the space of bilateral directed paths in $\Ac$.
Namely,
letting $A=(A(e, e'))_{e, e' \in E}$ be the adjacency matrix (for directed edges) of $\Ac$,
where $A(e, e')=1$ if $(e, e')$ is a directed path in $\Ac$, and $0$ otherwise,
we define
\[
\SS:=\Big\{(\o_i)_{i \in \Z} \ : \ A(\o_i, \o_{i+1})=1 \quad \text{for all $i \in \Z$}\Big\}.
\]
Let us define the shift $\s$ on $\SS$,
\[
\s:\SS \to \SS, \quad \s(\o_i)_{i \in \Z}:=(\o_{i+1})_{i\in \Z}.
\]
Then we obtain a subshift of finite type $(\SS, \s)$ with the set of alphabets $E$.
Let us define the metric in $\SS$ by
\[
d_\SS(\o, \o'):=e^{-n} \quad \text{where $n:=\sup\{k \ge 0 \ : \ \o_i=\o'_i \ \text{for all $|i| \le k$} \}$ and $e^{-\infty}=0$}.
\]

We consider the map assigning to a bilateral path $\o=(\o_i)_{i \in \Z}$ in $\SS$ a pair of points $\x_-(\o)$ and $\x_+(\o)$ in $\partial \G$,
where $\x_-(\o)$ and $\x_+(\o)$ are 
the extremes of geodesic rays 
\[
(o, w(\o_{-1})^{-1}, w(\o_{-1})^{-1}w(\o_{-2})^{-1}, \dots) \quad \text{and} \quad (o, w(\o_0), w(\o_0)w(\o_1), \dots),
\]
respectively.
Abusing the notation, we denote this map by
\[
w_\ast: \SS \to \partial^2 \G, \quad \o \mapsto (\x_-(\o), \x_+(\o)).
\]
The map $w_\ast$ is continuous; although it is not injective nor surjective, it is $C$-to-$1$ for some constant $C \ge 1$.
\begin{lemma}\label{Lem:Cto1}
Let $\rho(\x, \x')=e^{-(\x|\x')_o}$ for a word metric $d_S$ and $\rho_\times((\x_-, \x_+), (\x_-', \x_+'))=\max\{\rho(\x_-, \x_-'), \rho(\x_+, \x_+')\}$. 
\begin{itemize}
\item[(i)]
There exists a constant $L >0$ such that
\[
\rho_\times(w_\ast(\o), w_\ast(\o')) \le L d_\SS(\o, \o') \quad \text{for all $\o, \o' \in \SS$}.
\]
\item[(ii)]
There exists an integer $C \ge 1$ such that
for every $(\x_-, \x_+) \in \partial^2 \G$,
the number of pre-image $w_\ast^{-1}(\x_-, \x_+)$ is at most $C$.
\end{itemize}
\end{lemma}

\proof
For any $\o, \o'\in \SS$, let $w_\ast(\o)=(\x_-, \x_+)$ and $w_\ast(\o')=(\x_-', \x_+')$.
If $d_\SS(\o, \o')=e^{-n}$ for $n \ge 0$,
then $(\x_+|\x_+')_o \ge n-2 \d$ and $(\x_-|\x_-')_o \ge n-2 \d$ by the $\d$-hyperbolicity of the Cayley graph $\Cay(\G, S)$,
and this shows (i).

Note that the $\d$-hyperbolicity of $\Cay(\G, S)$ implies that two geodesics with the same extreme points are within $C_\d$-Hausdorff distance for some constant $C_\d \ge 0$.
For $(\x_-, \x_+) \in \partial^2 \G$,
if there exists an $\o \in \SS$ such that $w_\ast(\o)=(\x_-, \x_+)$,
then for any other $\o' \in \SS$ with $w_\ast(\o')=(\x_-, \x_+)$, 
the corresponding geodesic $w(\o')$ in $\Cay(\G, S)$ passes through $o$ from $\x_-$ to $\x_+$.
Since $\o'$ is a path in the automaton $\Ac$, the number of geodesic rays from $o$ of the form $w((\o')_{i=0, 1, \dots})$ following $w(\o)$ within $C_\d$-distance is at most $|B_S(o, C_\d)|$.
The same is true for geodesic rays from $o$ of the form $w((\o')_{i=-1, -2, \dots})$.
Therefore
the number of such $\o'$ is at most $|B_S(o, C_\d)|^2$, and this implies (ii).
\qed

Associated to the cocycle $\a: \G \times \partial^2 \G \to \R$,
we define a function on $\SS$ by
\[
\wt \a:\SS \to \R, \quad \o \mapsto \a(s_0^{-1}, w_\ast (\o)) \quad \text{where $s_0:=w(\o_0)$ for $\o=(\o_i)_{i \in \Z}$}.
\]

\begin{lemma}\label{Lem:suspension}
There exist an $N \ge 1$ and positive constants $c_1, c_2>0$ such that 
\[
S_N\wt \a(\o):=\sum_{n=0}^{N-1}\wt \a(\s^n\o) \in [c_1, c_2] \quad \text{for all $\o \in \SS$}.
\]
\end{lemma}

\proof
For any integer $N \ge 1$, the cocycle identity of $\a$ yields 
\begin{align*}
S_N\wt \a(\o)	&=\a(s_0^{-1}, w_\ast(\o))+\a(s_1^{-1}, s_0^{-1}w_\ast(\o))+\cdots+\a(s_{N-1}^{-1}, s_{N-2}^{-1}\cdots s_0^{-1}w_\ast(\o))\\
			&=\a((s_0\cdots s_{N-1})^{-1}, w_\ast(\o)),
\end{align*}
where $s_i=w(\o_i)$ for $i=0, \dots, N-1$.
Letting $(\x_-, \x_+)=w_\ast(\o)$,
we take a $C_\d$-rough geodesic $\phi_{-,+}$ from $\x_-$ to $\x_+$ with $\wh d(\phi_{-,+}(0), o) \le C_\d$. 
Then we have
\[
\a(\phi_{-,+}(r)^{-1}, \x_-, \x_+)=r\pm C_\d, \quad \text{for $r \in \R$}.
\]
Since for all large enough $N$, there exists an $r>2C_\d$ such that $\wh d(\phi_{-,+}(r), s_1\cdots s_N) \le C_\d$.
Hence 
$\a((s_0\cdots s_{N-1})^{-1}, w_\ast(\o)) \ge \a(\phi_{-,+}(r)^{-1}, \x_-, \x_+)-2C_\d$ and
there exist $N \ge 1$ and $c>0$ such that $S_N\wt \a(\o) \ge c$ for all $\o \in \SS$.
Letting 
\[
C:=\sup\{|\a(x, w_\ast \o)| \ : \ \o \in \SS, \ x \in B_S(o, N)\}<\infty,
\]
we obtain $S_N\wt \a(\o) \le C$ for all $\o \in \SS$.
\qed

\subsection{Coding the flow space}

Fix an integer $N \ge 1$ in Lemma \ref{Lem:suspension}.
Let 
\[
\Sus(\SS, S_N\wt \a):=\(\SS \times \R\)/\sim,
\]
where 
\[
(\o, t+S_N\wt \a(\o)) \sim (\s^N\o, t) \quad \text{for $(\o, t) \in \SS \times \R$}.
\]
We endow $\Sus(\SS, S_N \wt \a)$ with the quotient topology from $\SS \times \R$.
If we define the $\R$-action on $\Sus(\SS, S_N\wt \a)$ by $\s_t[\o, s]:=[\o, t+s]$ for $t \in \R$ and $[\o, s] \in \Sus(\SS, S_N\wt \a)$,
then $\Sus(\SS, S_N\wt \a)$ defines a {\it suspension flow} over $(\SS, \s^N)$.

Let us define $\wt \Pi:\SS\times \R \to \partial^2 \G \times \R$ by
\[
\wt \Pi(\o, t):=(w_\ast(\o), t) \quad \text{for $(\o, t) \in \SS \times \R$},
\]
and the map by
\[
\Pi: \Sus(\SS, S_N\wt \a) \to \Fc_\a, \quad [\o, t] \mapsto [w_\ast(\o), t].
\]
We shall show that this map $\Pi:\Sus(\SS, S_N\wt \a) \to \Fc_\a$ is well-defined.
For $(\o, t) \in \SS \times \R$,
we have
\[
\wt \Pi(\o, t+S_N\wt \a(\o))=(w_\ast(\o), t+S_N\wt \a(\o)) \quad \text{and} \quad \wt \Pi(\s^N\o, t)=(x_{N}^{-1}\cdot w_\ast(\o), t),
\]
where $x_N:=w(\o_0)\cdots w(\o_{N-1})$.
By definition of the $(\G, \a)$-action on $\partial^2\G \times \R$,
\begin{align*}
x_N^{-1}\cdot \wt \Pi(\o, t+S_N\wt \a(\o))	=(x_N^{-1}\cdot w_\ast(\o), t+S_N\wt \a(\o)-\a(x_N^{-1}, w_\ast(\o)))=\wt \Pi(\s^N\o, t)
\end{align*}
where we have used $S_N \wt \a(\o)=\a(x_N^{-1}, w_\ast(\o))$ in the last equality.
This shows that $\wt \Pi$ defines the map $\Pi:\Sus(\SS, S_N\wt \a) \to \Fc_\a$.

\begin{proposition}\label{Prop:coding}
The map
\[
\Pi: \Sus(\SS, S_N\wt \a) \to \Fc_\a, \quad [\o, t] \mapsto [w_\ast(\o), t]
\]
is continuous and satisfies that $\Pi\circ \s_t=\Phi_t \circ \Pi$ for $t \in \R$.
Moreover, $\Pi$ is a surjective map and there exists a constant $C\ge 1$ such that
the cardinality of $\Pi^{-1}(x)$ is at most $C$ for all $x \in \Fc_\a$.
\end{proposition}

\proof
By definition $\Pi[\o, t]=[\wt \Pi(\o, t)]$ for $(\o, t) \in \SS \times \R$, and since $\wt \Pi$ is continuous, $\Pi$ is a continuous map.
Furthermore,
if we define an $\R$-action on $\SS \times \R$ by $\wt \s_t(\o, s):=(\o, s+t)$,
then $\s_t[\o, s]=[\wt \s_t(\o, s)]$ and
\[
\wt \Pi \circ \wt \s_t=\wt \F_t \circ \wt \Pi \quad \text{for $t \in \R$}.
\]
This relation descends to $\Pi \circ \s_t=\F_t \circ \Pi$ for all $t \in \R$.

Let us show that the map $\Pi: \Sus(\SS, S_N\wt \a) \to \Fc_\a$ is surjective.
For given $(\x_-, \x_+) \in \partial^2 \G$,
let $\phi=(\phi(n))_{n \in \Z}$ be a geodesic from $\x_-$ to $\x_+$ in the Cayley graph $\Cay(\G, S)$ where $\phi(n) \in \G$ for $n \in \Z$.
For $n \ge 1$, let $\x_{-n}:=\phi(-n)$ and $\x_n:=\phi(n)$, for which we have $\x_{-n} \to \x_-$ and $\x_n \to \x_+$ in $\G \cup \partial \G$ as $n \to \infty$.
The strongly Markov automatic structure $(\Ac, w)$ gives a unique geodesic path $\phi_n$ from $\x_{-n}$ to $\x_n$ as the image of a directed path of length $2n+1$ in $\Ac$ from the initial state $s_\ast$.
We parametrize $\phi_n$ so that $\phi_n:[-n, n] \to \Cay(\G, S)$ with $\x_{-n}=\phi_n(-n)$ and $\x_n=\phi_n(n)$ and extend $\phi_n$ on $\Z$ in such a way that $\phi_n$ and $\phi$ coincide on $(-\infty, -n]\cup[n, \infty)$.
The $\d$-hyperbolicity of $\Cay(\G, S)$ implies that $d_S(\phi(i), \phi_n(i)) \le C_\d$ for all $i \in [-n, n]$.
Hence there exists a subsequence $\phi_{n_k}$ such that $\phi_{n_k}(i)$ converges to a $\g_i \in \G$ for each $i \in \Z$ as $n_k \to \infty$.
The path $(\g_i)_{i \in \Z}$ is a geodesic from $\x_-$ and $\x_+$, and since any finite subpath $(\g_i)_{-N \le i \le N}$ is given by the image of a directed path in $\Ac$,
there exists an $\o \in \SS$ such that $w(\o_i)=\g_{i}^{-1}\cdot \g_{i+1}$ for all $i \in \Z$.
Then, we have
\[
w_\ast(\o)=(\g_0^{-1}\x_-, \g_0^{-1}\x_+).
\]
This shows that for an arbitrary $(\x_-, \x_+, t) \in \partial^2 \G \times \R$, 
there exist an $\o \in \SS$ and a $\g_0 \in \G$ such that
$\g_0\cdot(w_\ast(\o), t+\a(\g_0, \g_0^{-1}\x_-, \g_0^{-1}\x_+))=(\x_-, \x_+, t)$.
Therefore
\[
\Pi [\o, t+\a(\g_0, \g_0^{-1}\x_-, \g_0^{-1}\x_+)]=[\x_-, \x_+, t],
\]
and the map $\Pi:\Sus(\SS, S_N\wt \a) \to \Fc_\a$ is surjective.

Let us show that the map $\Pi$ is $C$-to-one for some constant $C \ge 1$.
Since $\Pi$ is surjective,
for any $[\x_-, \x_+, s] \in \Fc_\a$, 
we have a $[\o, s] \in \Sus(\SS, S_N\wt \a)$ such that $\Pi[\o, s]=[\x_-, \x_+, s]$.
Letting
\[
Z:=\Big\{(\o, t) \in \SS \times \R \ : \ 0 \le t <S_N\wt \a(\o)\Big\},
\]
we identify $Z$ with $\Sus(\SS, S_N\wt \a)$.
Fix a unique lift $(\o, s) \in Z$ of $[\o, s]$.
We count the number of $(\o', s') \in Z$ such that $\wt \Pi(\o', s')$ and $\wt \Pi(\o, s)$ are in the same $(\G, \a)$-orbit in $\partial^2 \G \times \R$.
Note that $\wt \Pi(Z)=\{(w_\ast(\o), t) \in \partial^2 \G \times \R \ : \ (\o, t) \in Z\}$
is included in a compact set $K$ in $\partial^2 \G \times \R$
since $0<S_N\wt \a(\o)\le C$ by Lemma \ref{Lem:suspension} and $\lbr \x_-|\x_+\rbr_o \le C'$ where $w_\ast(\o)=(\x_-, \x_+)$ for all $\o \in \SS$.
If there exists a $\g \in \G$ such that $\g\cdot (w_\ast(\o), s)=(w_\ast(\o'), s')$, then the number of such $\g$ is at most $C_K$, which is a constant depending only on $K$ since the $(\G, \a)$-action on $\partial^2 \G \times \R$ is properly discontinuous by Lemma \ref{Lem:Gamma-action}.
For each such $\g \in \G$, we have
\[
\Big(w_\ast(\o'), s'\Big)=\Big(\g \cdot \x_-, \g \cdot \x_+, s+\a(\g, \x_-, \x_+)\Big),
\]
where $w_\ast(\o)=(\x_-, \x_+)$.
There exists a constant $C \ge 1$ such that the number of $\o'$ such that $w_\ast(\o')=(\g\cdot \x_-, \g \cdot \x_+)$ is at most $C$ by Lemma \ref{Lem:Cto1} (ii).
Therefore the number of $(\o', s') \in Z$ such that $\Pi[\o', s']=\Pi[\o, s]$ is at most $C_K \cdot C$.
We conclude that $\Pi:\Sus(\SS, S_N\wt \a) \to \Fc_\a$ is $(C_K\cdot C)$-to-one.
\qed

\subsection{Coding the measures}

Let $\Pi:\Sus(\SS, S_N \wt \a) \to \Fc_\a$ be the map in Proposition \ref{Prop:coding}.
For the simplicity of notation, fix $N$ and we write $r:=S_N\wt \a$ and $\Sus(\SS, r):=\Sus(\SS, S_N\wt \a)$.

\begin{lemma}\label{Lem:lift-measure}
For any $\{\F_t\}_{t \in \R}$-invariant Borel probability measure $m$ on $\Fc_\a$, there exists a $\{\s_t\}_{t \in \R}$-invariant Borel probability measure $\wt m$ on $\Sus(\SS, r)$ such that
\[
\Pi_\ast \wt m=m.
\]
\end{lemma}

\proof
Let $m$ be any probability measure on $\Fc_\a$ such that $m \circ \F_{-t}=m$ for all $t \in \R$.
If we have a continuous function $f \in C(\Fc_\a)$, then $f\circ \Pi \in C(\Sus(\SS, r))$, and
the positive linear functional
\[
L(f\circ \Pi):=\int_{\Fc_\a}f\,dm
\]
extends to a linear functional $\wt L$ on $C(\Sus(\SS, r))$ such that $|\wt L(\wt f)| \le \|\wt f\|_\infty$ for all $\wt f \in C(\Sus(\SS, r))$ by the Hahn-Banach theorem.
Note that $\wt L(1)=1$, and moreover $\wt L$ is positive. 
Indeed, suppose that there is an $\wt f$ on $\Sus(\SS, r)$ such that $\wt f \ge 0$ and $\wt L(\wt f)<0$, 
then normalizing $\wt f$ so that $0 \le \wt f \le 1$ if necessary, we have
\[
\wt L(1-\wt f)=1-\wt L(\wt f)>1;
\]
since $\|1-\wt f\|_\infty \le 1$, this is a contradiction.
Hence by the Riesz representation theorem there exists a Borel probability measure $\wt \m$ on $\Sus(\SS, r)$ such that
\[
\int_{\Sus(\SS, r)}f \circ \Pi\,d\wt \m=\int_{\Fc_\a}f\,dm \quad \text{for $f \in C(\Fc_\a)$}.
\]
For all $T>0$, let
$\wt \m_T:=\frac{1}{2T}\int_{-T}^T(\s_{t})_\ast \wt \m\,dt$.
Then there exists a sequence $T_i \to \infty$ such that $\wt \m_{T_i}$ weak-star converges to a probability measure $\wt m$, which is $\{\s_t\}_{t \in \R}$-invariant.
Since $\Pi \circ \s_t=\F_t\circ \Pi$ for all $t \in \R$, 
we have
$\Pi_\ast \wt m=m$ as required.
\qed

For any real values $a<b$,
we denote by $\Leb_{[a, b)}$ the Lebesgue measure restricted on the interval $[a, b)$ in $\R$.
For a probability measure $\lambda$ on $\SS$, let $w_\ast \lambda:=\lambda \circ w_\ast^{-1}$ be the pushforward of $\lambda$ by $w_\ast$ on $\partial^2 \G$.

\begin{lemma}\label{Lem:Nto1}
For any $\s^N$-invariant Borel probability measure $\lambda$ on $\SS$, and for any $T_0>0$,
there exist a $\s$-invariant Borel probability measure $\wt \lambda$ on $\SS$ and some $T \ge T_0$
such that
the following inequalities hold on $\partial^2 \G \times \R$, 
\[
\frac{1}{N}w_\ast \lambda \otimes \Leb_{[0, T_0)} \le w_\ast \wt \lambda \otimes \Leb_{[0, T_0)}
\le \frac{1}{N}\sum_{|\g|_S \le N-1}\g.\(w_\ast \lambda \otimes \Leb_{[-T, T)}\).
\]
\end{lemma}

\proof
Let $\wt \lambda:=\frac{1}{N}\sum_{i=0}^{N-1}\s^i_\ast \lambda$.
Then $\wt \lambda$ is a $\s$-invariant Borel probability measure on $\SS$.
Since $\lambda \le N\cdot \wt \lambda$, we have
$w_\ast \lambda \otimes \Leb_{[0, T_0)} \le N w_\ast \wt \lambda \otimes \Leb_{[0, T_0)}$.

On the other hand, we note that for $\o \in \SS$, 
\[
w_\ast(\s \o)=(\x_-(\s \o), \x_+(\s \o))=s_0^{-1}.(\x_-(\o), \x_+(\o))=s_0^{-1}.w_\ast(\o),
\]
where $s_0=w(\o_0)$.
For any Borel set $A$ in $\partial^2 \G$ and for each $s \in S$,
we have
\[
\lambda\(\s^{-1}w_\ast^{-1}A \cap\{w(\o_0)=s\}\)
\le \lambda\(\{s^{-1}.w_\ast (\o) \in A\}\)
=w_\ast \lambda(sA)
=(s^{-1}). w_\ast \lambda(A),
\]
and thus since $S=S^{-1}$,
$w_\ast \s_\ast \lambda \le \sum_{s \in S}s.w_\ast \lambda$.
In fact, we have for all integer $i \ge 0$,
$w_\ast \s^i_\ast \lambda \le \sum_{|\g|_S =i}\g.w_\ast \lambda$ and
\[
w_\ast \wt \lambda \le \frac{1}{N}\sum_{|\g|_S \le N-1}\g.w_\ast \lambda.
\]
Let $T:=T_0+T'$ where
\[
T':=\sup\Big\{|\a(\g, \x_-, \x_+)| \ : \ |\g|_S \le N-1, \ (\x_-|\x_+)_o \le 4\d\Big\}<\infty.
\]
For any compactly supported continuous function $f \ge 0$ on $\partial^2 \G \times \R$,
we have
\begin{align*}
\int_{\partial^2 \G \times \R}f(\x, t)\,d\((\g.w_\ast \lambda)\otimes \Leb_{[0, T_0)}\)
&=\int_{\partial^2 \G \times \R}f(\g.\x, t)\,d\(w_\ast \lambda \otimes \Leb_{[0, T_0)}\)\\
&=\int_{\partial^2 \G \times \R}f\circ \g(\x, t-\a(\g, \x))\,d\(w_\ast \lambda \otimes \Leb_{[0, T_0)}\)\\
&\le \int_{\partial^2 \G \times \R}f\circ \g(\x, t)\,d\(w_\ast \lambda \otimes \Leb_{[-T, T)}\),
\end{align*}
for all $\g \in \G$ such that $|\g|_S \le N-1$,
where we have used the definition of $T=T_0+T'$ in the last inequality.
Therefore we obtain
\[
w_\ast \wt \lambda \otimes \Leb_{[0, T_0)} \le \frac{1}{N}\sum_{|\g|_S \le N-1}\g.\(w_\ast \lambda \otimes \Leb_{[-T, T)}\),
\]
and thus we conclude the claim.
\qed

Let us consider
\[
Z:=\Big\{(\o, t) \in \SS \times \R \ : \ 0 \le t <r(\o)\Big\},
\]
and identify $\Sus(\SS, r)$ with $Z$ and also the measures on them.
Let $\Mcc(\{\s_t\}_{t \in \R})$ be the set of flow $\{\s_t\}_{t\in \R}$-invariant Borel probability measures on $\Sus(\SS, r)$.
There is a bijection between $\Mcc(\s^N, \SS)$ and $\Mcc(\{\s_t\}_{t \in \R})$,
\[
\Mcc(\s^N, \SS) \to \Mcc(\{\s_t\}_{t \in \R}), \quad \lambda \mapsto \n_\lambda:=\frac{1}{\lambda \otimes \Leb(Z)}\lambda\otimes \Leb|_Z,
\]
where $\Leb$ stands for the Lebesgue measure on $\R$.
The reverse map is given by for $\n \in \Mcc(\{\s_t\}_{t \in \R})$ taking $\wbar \n$ on $\SS$ by disintegration $\n=\int_\SS \Leb_{[0, r(\o))}\,d\wbar \n$
and assigning $\wbar \n(\SS)^{-1}\cdot \wbar \n$.

\begin{proposition}\label{Prop:coding_measure}
There exist constants $T_0, c_1, c_2>0$ such that the following hold.
\begin{itemize}
\item[(i)] 
For every flow $\{\F_t\}_{t \in \R}$-invariant Borel probability measure $m$ on $\Fc_\a$,
there exists a $\s$-invariant Borel probability measure $\wt \lambda$ on $\SS$ satisfying that
\[
c_1 m \le \pi_\ast \(w_\ast \wt \lambda \otimes \Leb_{[0, T_0)}\) \le c_2 m,
\]
where $\pi:\partial^2 \G \times \R \to \Fc_\a$ is the quotient map.
\item[(ii)]
For every $\G$-invariant Radon measure $\L$ on $\partial^2 \G$,
there exists a $\s$-invariant Borel probability measure $\wt \lambda$ on $\SS$ satisfying that
\[
c_1 \L \otimes dt \le \sum_{\g \in \G}\g.\(w_\ast \wt \lambda \otimes \Leb_{[0, T_0)}\) \le c_2 \L \otimes dt \quad \text{on $\partial^2 \G \times \R$}.
\]
\end{itemize}
\end{proposition}

Note that for any compactly supported finite measure $\m$ on $\partial^2 \G \times \R$,
the sum $\sum_{\g \in \G}\g.\m$ is well-defined since the $(\G, \a)$-action on $\partial^2 \G \times \R$ is properly discontinuous.

\def\pr{{\rm pr}}

\proof
First we show (ii) by using (i).
Let $\L$ be a $\G$-invariant Radon measure on $\partial^2 \G$.
Letting $\i:\Fc_\a \to D$ be a Borel section,
we have a $\{\F_t\}_{t \in \R}$-invariant finite measure $m$ on $\Fc_\a$ such that
\[
\L \otimes dt =\sum_{\g \in \G}\g.\i_\ast m,
\]
by Lemma \ref{Lem:quotient-measure} (and Remark \ref{Rem:quotient-measure}).
By (i), there exists a $\s$-invariant probability measure $\wt \lambda$ on $\SS$ satisfying
that $c_1 m \le \pi_\ast\(w_\ast \wt \lambda \otimes \Leb_{[0, T_0)}\) \le c_2 m$.
Since we have that
\[
\sum_{\g \in \G}\g. \i_\ast \pi_\ast\(w_\ast \wt \lambda \otimes \Leb_{[0, T_0)}\)=\sum_{\g \in \G}\g. \(w_\ast \wt \lambda \otimes \Leb_{[0, T_0)}\),
\] 
we obtain (ii),
\[
c_1 \L \otimes dt \le \sum_{\g \in \G}\g. \(w_\ast \wt \lambda \otimes \Leb_{[0, T_0)}\) \le c_2 \L \otimes dt.
\]

Next let us show (i).
Let $m$ be a $\{\F_t\}_{t \in \R}$-invariant probability measure on $\Fc_\a$.
First we shall show that the claim holds for a $\s^N$-invariant probability measure $\lambda$ on $\SS$ for an integer $N \ge 1$.
By Lemma \ref{Lem:lift-measure}, there exists a $\{\s_t\}_{t \in \R}$-invariant Borel probability measure $\wt m$ on $\Sus(\SS, r)$
such that $\Pi_\ast \wt m=m$,
where $r=S_N\wt \a$ for some $N \ge 1$.
Identifying $\Sus(\SS, r)$ with $Z$,
we write the measure $\wt m$ by $\frac{1}{\lambda\otimes \Leb (Z)}(\lambda \otimes \Leb)|_Z$ where $\lambda$ is a $\s^N$-invariant Borel probability measure on $\SS$.
Formally, letting $\pr: \SS \times \R \to \Sus(\SS, r)$ be the quotient map,
we have
\[
\wt m=\pr_\ast \frac{1}{\lambda\otimes \Leb (Z)}(\lambda \otimes \Leb)|_Z.
\]
Let $T_0:=\max_{\o \in \SS} r(\o)$ and $c_0:=\min_{\o \in \SS} r(\o)$, where we have $c_0>0$.
Note that $c_0 \le (\lambda \otimes \Leb)(Z) \le T_0$.
By the definition of $\wt \Pi$, it holds that $\pi\circ \wt \Pi=\Pi \circ \pr$.
Therefore we have
\[
m=\Pi_\ast \pr_\ast \frac{1}{\lambda \otimes \Leb(Z)}(\lambda \otimes \Leb)|_Z
=\pi_\ast \wt \Pi_\ast \frac{1}{\lambda \otimes \Leb(Z)}(\lambda \otimes \Leb)|_Z.
\]
By using disintegration with respect to $w_\ast \lambda$, we obtain
\[
\wt \Pi_\ast (\lambda \otimes \Leb)|_Z=\int_{\partial^2 \G} l_\x\, d(w_\ast \lambda) (\x),
\]
where $l_\x=\sum_{w_\ast\o=\x}\Leb_{[0, r(\o))}$ for $w_\ast \lambda$-almost every $\x$ in $\partial^2 \G$.
Noting that the map $w_\ast: \SS \to \partial^2 \G$ is $C$-to-one for some $C \ge 1$ by Lemma \ref{Lem:Cto1} (ii),
we have
\[
\wt \Pi_\ast(\lambda \otimes \Leb)|_Z \le C\cdot w_\ast \lambda \otimes \Leb_{[0, T_0)}.
\]
This shows that 
\begin{equation}\label{Eq:coding_measure1}
m \le \frac{C}{c_0}\pi_\ast (w_\ast \lambda \otimes \Leb_{[0, T_0)}).
\end{equation}
We shall show the other estimate.
Since $\wt m$ is $\{\s_t\}_{t \in \R}$-invariant, for any $T \ge 0$, $(T+c_0)\cdot \wt m=\int_{-c_0}^T (\s_t)_\ast \wt m$.
Noting that $(\lambda \otimes \Leb)|_Z \ge \lambda \otimes \Leb_{[0, c_0)}$, we have that
for any $T \ge 0$,
\[
\int_{-c_0}^T(\s_t)_\ast (\lambda \otimes \Leb)|_Z \ge c_0\cdot \lambda \otimes \Leb_{[0, T)}
\quad
\text{and}
\quad
\wt m \ge \frac{c_0}{T_0(T+c_0)}\pr_\ast (\lambda \otimes \Leb_{[0, T)}).
\]
Hence by using the relation $\pi\circ \wt \Pi=\Pi \circ \pr$, we obtain for all $T\ge 0$,
\begin{align}\label{Eq:coding_measure2}
m=\Pi_\ast \wt m 
&\ge \frac{c_0}{T_0(T+c_0)} \Pi_\ast \pr_\ast (\lambda \otimes \Leb_{[0, T)}) \nonumber\\
&= \frac{c_0}{T_0(T+c_0)} \pi_\ast \wt \Pi_\ast(\lambda \otimes \Leb_{[0, T)}) 
\ge \frac{c_0}{T_0(T+c_0)} \pi_\ast (w_\ast \lambda \otimes \Leb_{[0, T)}).
\end{align}

Finally,
applying Lemma \ref{Lem:Nto1}, we have a $\s$-invariant probability measure $\wt \lambda$ on $\SS$
such that for some $T \ge T_0$,
\begin{equation}\label{Eq:coding_measure3}
\frac{1}{N}w_\ast \lambda \otimes \Leb_{[0, T_0)} \le w_\ast \wt \lambda \otimes \Leb_{[0, T_0)}
\le \frac{1}{N}\sum_{|\g|_S \le N-1}\g.\(w_\ast \lambda \otimes \Leb_{[-T, T)}\).
\end{equation}
By \eqref{Eq:coding_measure1} and the first inequality of \eqref{Eq:coding_measure3}, we have
$c_1 m \le \pi_\ast\(w_\ast \wt \lambda \otimes \Leb_{[0, T_0)}\)$.
By \eqref{Eq:coding_measure2}, we have that for some positive constant $c>0$,
\[
3Tm=\int_{-2T}^T(\F_t)_\ast m\,dt \ge c \int_{-2T}^T(\F_t)_\ast\pi_\ast(w_\ast \lambda \otimes \Leb_{[0, T)})\,dt 
\ge cT\pi_\ast\(w_\ast\lambda \otimes \Leb_{[-T, T)}\).
\]
Since $\pi\circ \g=\pi$ for all $\g \in \G$, by the second inequality of \eqref{Eq:coding_measure3} we obtain
\begin{align*}
\pi_\ast\(w_\ast \wt \lambda \otimes \Leb_{[0, T_0)}\)
\le \frac{|B_S(o, N-1)|}{N}\pi_\ast\(w_\ast \lambda \otimes \Leb_{[-T, T)}\)
\le \frac{|B_S(o, N-1)|}{N}\cdot \frac{3}{c}m,
\end{align*}
where $B_S(o, N-1)$ is the ball of radius $N-1$ centered at $o$ in $\Cay(\G, S)$.
We conclude the claim.
\qed

\section{Symbolic dynamics}\label{Sec:symbolic}

Recall that
associated to the cocycle $\a: \G \times \partial^2 \G \to \R$,
we have defined the function
\[
\wt \a:\SS \to \R, \quad \o \mapsto \a(s_0^{-1}, w_\ast (\o)),
\]
where $s_0:=w(\o_0)$ for $\o=(\o_i)_{i \in \Z}$.
Similarly, we define
\[
\wt b: \SS \to \R, \quad \o \mapsto \wh b_o(s_0, \x_+(\o)) 
\quad \text{and} \quad
u : \SS \to \R, \quad \o \mapsto \lbr \x_-(\o) | \x_+(\o) \rbr_o.
\]

\begin{lemma}\label{Lem:alpha-Holder}
The functions $\wt \a$, $\wt b$ and $u$ are H\"older continuous on $(\SS, d_\SS)$, i.e.,
there exist constants $L \ge 0$ and $s>0$ such that
\[
|\wt \a(\o)-\wt \a(\o')| \le L\, d_\SS(\o, \o')^s \quad \text{for all $\o, \o' \in \SS$},
\]
and the same estimates hold for $\wt b$ and $u$.
Moreover, we have
\[
\wt \a(\o)=-\wt b(\o)+u(\s \o)-u(\o) \quad \text{for all $\o \in \SS$}.
\]
\end{lemma}

\proof
For any $n \ge 0$,
let $\o, \o' \in \SS$ be any two sequences satisfying that $\o \neq \o'$ and $d_\SS(\o, \o')=e^{-n}$.
Since the Busemann function $\wh b_o(x, x_n) \to \wh b_o(x, \x)$ as $x_n \to \x$ in $\G \cup \partial \G$ for each $x \in \G$ and $\wh d$ is strongly hyperbolic \eqref{Eq:strongly_hyperbolic_word} in Section \ref{Sec:strongly_hyperbolic}, there exist constant $L \ge 0$ and $s>0$ such that
\begin{align*}
&|\wh b_o(s_0, \x_+(\o))-\wh b_o(s_0, \x_+(\o'))| \le Le^{-s(n-1)},
\end{align*}
where $s_0=w(\o_0)$, and thus
$|\wt b(\o)-\wt b(\o')| \le Le^s e^{-n}=Le^s\,d_\SS(\o, \o')$.
Moreover, the strong hyperbolicity of $\wh d$ implies that
\[
|\wh d(x, y)-\wh d(x', y)-\wh d(x, o)+\wh d(x', o)| \le Le^{-s n}
\]
if two geodesic segments connecting $x$ and $y$, and $x'$ and $o$, respectively, have a common geodesic of length $n$ in their $C_\d$-neighborhoods in the Cayley graph $\Cay(\G, S)$.
Adding $\wh d(y, o)$ and $-\wh d(y, o)$, we have that
\[
|2\lbr \x_-(\o')|\x_+(\o)\rbr_o - 2\lbr \x_-(\o)|\x_+(\o)\rbr_o| \le Le^{-s n}
\]
and 
\[
|2\lbr \x_-(\o')|\x_+(\o')\rbr_o - 2\lbr \x_-(\o')|\x_+(\o)\rbr_o| \le Le^{-s n}.
\]
This shows that
\[
|u(\o)-u(\o')| \le Le^{-sn}=L\,d_\SS(\o, \o').
\]
By the definition of $\a$ and the relation between $\wh b_0$ and the Gromov product, we have
\[
\wt \a(\o)=\a(s_0^{-1}, \x_-(\o), \x_+(\o))=\frac{1}{2}(\wh b_o(s_0, \x_-(\o))-\wh b_o(s_0, \x_+(\o)))
\]
and
\[
\wh b_o(s_0, \x_-(\o))+\wh b_o(s_0, \x_+(\o))=2\lbr \x_-(\o) | \x_+(\o)\rbr_{s_0}-2\lbr \x_-(\o) | \x_+(\o)\rbr_o.
\]
and thus
$\wt \a(\o)=-\wt b(\o)+u(\s \o)-u(\o)$.
Therefore $\wt \a$ is also H\"older continuous.
\qed

\subsection{Thermodynamic formalism}\label{Subsec:thermodynamic}
In this section, we work on the set of unilateral paths $\wbar \SS^+$.
It is better suited when we construct a shift-invariant measure on the set of bilateral paths $\SS$ by using transfer operators.
We introduce thermodynamic formalism on $\wbar \SS^+$ following \cite[Section 3.3]{GouezelLocalLimit}.

Let $\s:\wbar \SS^+ \to \wbar \SS^+$ be the shift defined by deleting the first edge of paths.
Define the metric in $\wbar \SS^+$ by
$d_{\wbar \SS^+}(\o, \o'):=e^{-n}$ where $n:=\sup\{k \ge 0 \ : \ \o_i=\o_i', \  0 \le i \le k\}$.
For every real-valued H\"older continuous function $\p:\wbar \SS^+ \to \R$ with respect to $d_{\wbar \SS^+}$,
we define the transfer operator $\Lc_\p$ acting on the space of continuous functions $f$ on $\wbar \SS^+$ by
\[
\Lc_\p f(\o)=\sum_{\s(\o')=\o}e^{\p(\o')}f(\o'),
\]
where the preimages of $\o=\emptyset$ by the shift $\s$ are paths of length $1$.
We will analyze the asymptotics of 
\[
\Lc_\p^n 1(\emptyset)=\sum_{\o' \,\text{of length $n$}}e^{S_n\p(\o')}, \quad \text{where $S_n \p(\o):=\p(\o)+\p(\s \o)+\cdots+\p(\s^{n-1}\o)$}.
\]
We call a finite directed graph {\it recurrent} when any vertex is accessible from any other vertex by a directed path, and {\it topologically mixing} when it is recurrent and there exists $N$ such that for all $n \ge N$ any two vertices are connected by a directed path of length $n$.
If the underlying graph of $\Ac$ is topologically mixing,
then the Ruelle-Perron-Frobenius theorem is applicable to describe the spectra of $\Lc_\p$ (\cite[Theorem 1.7]{Bowen}, \cite[Theorem 2.2]{ParryPollicott} and see also \cite[Theorem 3.6]{GouezelLocalLimit}).
If the graph $\Ac$ is only recurrent, then there is a period $p \ge 1$; every loop has the length of a multiple of $p$.
In the case when $p>1$, the set of vertices $V$ of $\Ac$ is a disjoint union of subsets $V_j$ for $j \in \Z/p\Z$ such that any edge of the origin in $V_j$ has the terminus in $V_{j+1}$.
We call such a decomposition $V=\bigsqcup_{j \in \Z/p\Z}V_j$ a {\it cyclic decomposition}.
Let $\wbar \SS_{j}^+$ be the set of path starting at a vertex in $V_j$ and the empty path.
Then $\s: \wbar \SS_{j}^+ \to \wbar \SS_{j+1}^+$ and,
the restriction of $\s^p$ on $\wbar \SS_j^+$ is a topological mixing subshift of finite type for each $j \in \Z/p\Z$.
%(In this case, $(\wbar \SS_j^+, \s^p)$ is regarded as a shift space based on alphabets consisting of a set of paths from $V_j$ to $V_j$ itself in $\Ac$.)
In the general case when $\Ac$ is not necessarily recurrent, we decompose $\Ac$ into components.
A {\it component} of $\Ac$ is a maximal subgraph which is recurrent.
For a general hyperbolic group, a strongly Markov automatic structure $\Ac$ can have several distinct components whose cardinality is greater than $1$.
(There are examples such as free groups and surface groups with the standard sets of generators admitting an $\Ac$ with a unique component which is not a singleton.)

To each component $\Cc$, 
we define the transfer operator $\Lc_\Cc$ by restricting $\p$ on the set of paths staying in $\Cc$.
Then $\Lc_\Cc$ has finitely many eigenvalues of maximal modulus ${\mathcal R}_\Cc$ and they are all simple and isolated.
Let $\Pr_\Cc(\p, \s):=\log {\mathcal R}_\Cc$.
(This value will be called the {\it pressure}.)
We define
\[
\Pr(\p, \s):=\max_\Cc \Pr_\Cc(\p, \s),
\]
where the maximum is taken over all components $\Cc$ of $\Ac$.
We call a component $\Cc$ {\it maximal} if $\Pr(\p, \s)=\Pr_\Cc(\p, \s)$.
A potential $\p$ is called {\it semisimple} if there is no directed path between any two different maximal components.

Let $1_{[E_\ast]}$ be the indicator function on $\wbar \SS^+$ taking the value $1$ on the set of paths starting at $s_\ast$ and $0$ otherwise.
Then we have
\[
\Lc_\p^n 1_{[E_\ast]}(\emptyset)=\sum e^{S_n\p(\o)},
\]
where the summation runs over all paths $\o$ of length $n$ starting at $s_\ast$.
The following lemma is used to compute $\Pr(\p, \s)$ and decide whether the potential $\p$ is semisimple.

\begin{lemma}\label{Lem:Lc}
Let $k \ge 1$ be an integer.
If there exists a path from $s_\ast$ in $\Ac$ passing through successively $k$ different maximal components,
then there exists a constant $C>0$ such that for all integers $n \ge 1$,
we have
\[
\Lc_\p^n 1_{[E_\ast]}(\emptyset)\ge Cn^{k-1}e^{n \Pr(\p, \s)}.
\]
Let $L\ge 1$ be an integer.
If there are $L$ components in a finite directed graph of $\Ac$, 
then there exists a constant $C>0$ such that for any integers $n \ge 1$,
\[
\Lc_\p^n 1_{[E_\ast]}(\emptyset) \le C n^L e^{n \Pr(\p, \s)}.
\]
\end{lemma}

\proof
The first claim is a special case of \cite[Lemma 3.7]{GouezelLocalLimit} and the proof of the second claim is in \cite[Lemma 4.7]{THaus}.
\qed

\def\Hol{{\rm H}}
Although we do not use it directly in this paper, it is instructive to record the following theorem by Gou\"ezel in order to understand the situation when the potential $\p$ is semisimple.
We denote by $\|\cdot\|_{\Hol}$ the H\"older norm with some fixed exponent (whose explicit value is not used).

\begin{theorem}[Theorem 3.8 in \cite{GouezelLocalLimit}]\label{Thm:Gouezel}
If $\p$ is a semisimple potential, and
$\Cc_1, \dots, \Cc_I$ are the maximal components with period $p_i$,
where $\Cc_i=\bigsqcup_{j \in \Z/p_i\Z}\Cc_{i, j}$ is a cyclic decomposition for each $i=1, \dots, I$,
then there exist H\"older continuous functions $h_{i, j}$ and measures $\lambda_{i, j}$ such that
$\int_{\wbar \SS^+} h_{i, j}\,d\lambda_{i, j}=1$ and for every H\"older continuous function $f$, and for all $n \ge 0$,
\[
\Big\|\Lc_\p^n f-e^{n\Pr(\p, \s)}\sum_{i=1}^I\sum_{j=0}^{p_i-1}\(\int_{\wbar \SS^+} f d\lambda_{i, (j-n \ {\rm mod} \ p_i)}\)h_{i, j}\Big\|_{\Hol} \le C\|f\|_{\Hol} \cdot e^{-n \e_0}e^{n \Pr(\p, \s)},
\]
where $C \ge 0$ and $\e_0>0$ are constants independent of $n$ and $f$.
For each $i=1, \dots, I$, 
\[
\m_i:=\frac{1}{p_i}\sum_{j=0}^{p_i-1}h_{i, j}\lambda_{i, j}
\]
is a $\s$-invariant ergodic probability measure, and supported on $\SS_i^+$.
\end{theorem}

Defining the Busemann function $\wh b_o$ on the whole space $\G \cup \partial \G$ by
\[
\wh b_o(x, z):=\wh d(x, z)-\wh d(o, z) \quad \text{for $x, z \in \G$},
\]
we let $\wt b:\wbar \SS^+ \to \R$ where
$\wt b(\o):=\wh b_o(w(\o_0), w_\ast \o)$.
Note that $\wt b$ is H\"older continuous on $\wbar \SS^+$ by Lemma \ref{Lem:alpha-Holder}.
We denote the exponential volume growth rate $\wh v:=\gr(\wh d\,)$ relative to $\wh d$.

\begin{lemma}\label{Lem:potential}
If we define the potential by $\p:=\wh v\cdot \wt b$,
then $\Pr(\p, \s)=0$.
\end{lemma}

\proof
For any integer $n \ge 1$ and for all path $\o \in \wbar \SS^+$ of length $n$, we have
\[
S_k \p(\o)=\wh v\cdot \wh b_o(s_0 \cdots s_{k-1}, w_\ast \o) \quad \text{for all $1 \le k \le n$},
\]
where $s_k=w(\o_k)$ for $k=0, \dots, n-1$, and thus
$S_n \p(\o)=-\wh v\cdot \wh d(o, w_\ast \o)$.
Therefore we have
\[
\Lc_\p^n 1_{[E_\ast]}(\emptyset)=\sum_{\o:\ \text{length $n$ from $s_\ast$}} e^{S_n \p(\o)}
=\sum_{x \in \bS_n}e^{-\wh v \cdot \wh d(o, x)},
\]
where $\bS_n=\{x \in \G \ : \ |x|_S=n\}$.
Since the last term is bounded from above and from below by some constants independent of $n \ge 1$ 
by Lemma \ref{Lem:regular-growth},
Lemma \ref{Lem:Lc} shows that $\Pr(\p, \s)=0$.
\qed

\begin{remark}
In fact, the first part of Lemma \ref{Lem:Lc} shows that the potential $\p=\wh v\cdot \wt b$ is semisimple, which we do not exploit in this paper.
\end{remark}

The following lemma is specific to word metrics.
Let $v_S:=\gr(d_S)$ be the exponential volume growth rate with respect to $d_S$.

\begin{lemma}\label{Lem:potential_word}
Let $\wt b_S(\o)\equiv -1$ for all $\o \in \wbar \SS^+$ and $\p_S:=v_S\cdot \wt b_S$.
Then $\Pr(\p_S, \s)=0$.
\end{lemma}

\proof
Note that $S_k \p_S(\o)=-v_S\cdot k$ for all $1 \le k \le n$ for all $\o \in \wbar \SS^+$ of length $n$.
The proof follows as in Lemma \ref{Lem:potential}.
\qed

\subsection{Variational principle}

Let $\p$ be a H\"older continuous function on $\SS$ where we assume that $\p$ depends only on coordinates of non-negative indices, i.e., $\p(\o)=\p(\o_0, \o_1, \dots)$.
For each component $\Cc$ of $\Ac$, 
let $\SS_{\Cc}$ be the set of bilateral paths all the time staying in $\Cc$.
Then it holds that
\[
\Pr_\Cc(\p, \s)=\lim_{n \to \infty}\frac{1}{n}\log \sum_{[\o_0, \dots, \o_{n-1}]}\exp\Big(S_{\o_0, \dots, \o_{n-1}} \p\Big),
\]
where $S_{\o_0, \dots, \o_{n-1}}\p:=\sup\{S_n \p(\o) \ : \ \o \in [\o_0, \dots, \o_{n-1}]\}$ for cylinder sets $[\o_0, \dots, \o_{n-1}]$ in $\SS_\Cc$.

For any subshift of finite type $(\SS, \s)$,
we denote by $\Mcc(\s, \SS)$ the set of all $\s$-invariant probability measures on $\SS$.
For any $\lambda \in \Mcc(\s, \SS)$, let $h(\s, \lambda)$ be the measure theoretical entropy of $(\SS, \s, \lambda)$
(see Section \ref{Subsec:entropy} for the definition).

\begin{proposition}\label{Prop:variational_principle}
Let $\p$ be a H\"older continuous potential on $\SS$.
Then for each component $\Cc$ of $\Ac$, we have that
\begin{equation}\label{Eq:variational_principle}
\Pr_{\Cc}(\p, \s)=\sup_{\lambda \in \Mcc(\s, \SS_{\Cc})}\Big\{h(\lambda, \s)+\int_{\SS_{\Cc}}\p\,d\lambda \Big\},
\end{equation}
and there exists a unique $\s$-invariant Borel probability measure $\m_\Cc$ on $\SS_\Cc$ for which the supremum is attained.
Moreover, there exist constants $c_1, c_2>0$ such that 
\begin{equation}\label{Eq:Gibbs_i}
\frac{\m_\Cc[\o_0, \dots, \o_{n-1}]}{\exp\(-n\Pr_\Cc(\p, \s)+S_n \p(\o)\)} \in [c_1, c_2],
\end{equation}
for all $\o \in [\o_0, \dots, \o_{n-1}]$ and for all $n \ge 1$.
\end{proposition}

\proof
For each component $\Cc$, let $\Cc=\bigsqcup_{j \in \Z/p_\Cc \Z}\Cc_{j}$ be a cyclic decomposition where $p_\Cc$ is the period of $\Cc$.
Let $\SS_{\Cc, j}$ be the space of bilateral paths $(\o_i)_{i \in \Z}$ staying in $\Cc$ such that the edge $\o_0$ starts in $\Cc_j$.
Then we have that $\SS_{\Cc}=\bigsqcup_{j=0}^{p_\Cc-1}\SS_{\Cc, j}$ and $\s: \SS_{\Cc, j} \to \SS_{\Cc, j+1}$ for each $j \in \Z/p_\Cc\Z$.
The (two-sided) shift space $(\SS_{\Cc, j}, \s^{p_\Cc})$ is a topologically mixing subshift of finite type.
Since $\p$ is H\"older continuous on $\SS$, $S_{p_\Cc}\p$ is H\"older continuous on $\SS_{\Cc, 0}$ relative to the metric $d_{\SS}$ restricted to $\SS_{\Cc, 0}$.
Therefore we have that
\begin{equation}\label{Eq:Gibbs_ij}
\Pr_{\Cc_{0}}(S_{p_\Cc}\p, \s^{p_\Cc})=\sup_{\lambda \in \Mcc(\s^{p_\Cc}, \SS_{\Cc, 0})}\Big\{h(\lambda, \s^{p_\Cc})+\int_{\SS_{\Cc, 0}}S_{p_\Cc}\p\,d\lambda \Big\},
\end{equation}
and the supremum in \eqref{Eq:Gibbs_ij} is attained by a measure 
if and only if it is the unique Gibbs measure $\m_{\Cc, 0}$ with potential $S_{p_\Cc}\p$ \cite[p.19, Theorem 1.22]{Bowen}.
In this setting, the measure $\m_{\Cc, 0}$ is 
a unique $\s^{p_\Cc}$-invariant probability measure $\m_{\Cc, 0}$ on $\SS_{\Cc, 0}$
for which there exist constants $c_1, c_2>0$
such that
for any integer $N \ge 1$ and any cylinder set 
$[\o_0, \dots, \o_{N p_\Cc-1}]$ in $\SS_{\Cc}$ of a path $(\o_0, \dots, \o_{N p_\Cc-1})$ starting in $\Cc_0$,
one has
\begin{equation}\label{Eq:Gibbs-ij}
\frac{\m_{\Cc, 0}[\o_0, \dots, \o_{N p_\Cc-1}] }{\exp{\(-N p_\Cc \cdot \Pr_{\Cc}(\p, \s)+S_{N p_\Cc}\p(\o)\)}} \in [c_1, c_2],
\end{equation}
for any $\o \in [\o_0, \dots, \o_{N p_\Cc-1}]$.

Note that there is a one-to-one correspondence between $\Mcc(\s, \SS_{\Cc})$ and $\Mcc(\s^{p_\Cc}, \SS_{\Cc, 0})$.
For $\lambda \in \Mcc(\s, \SS_{\Cc})$, since $\lambda(\SS_{\Cc, j})=1/p_\Cc$ for any $j \in \Z/p_\Cc\Z$, letting the restriction $\lambda_0:=(p_\Cc \lambda)|_{\SS_{\Cc, 0}}$ 
we obtain $\lambda_0 \in \Mcc(\s^{p_\Cc}, \SS_{\Cc, 0})$.
On the other hand, for $\lambda_0 \in \Mcc(\s^{p_\Cc}, \SS_{\Cc, 0})$, letting $\lambda:=(1/p_\Cc)\sum_{j=0}^{p_\Cc-1}\lambda_0 \circ \s^{-j}$,
we obtain $\lambda \in \Mcc(\s, \SS_{\Cc})$.
These processes are mutually invertible.
Moreover, we have $\Pr_{\Cc_{0}}(S_{p_\Cc}\p, \s^{p_\Cc})=p_\Cc\, \Pr_{\Cc}(\p, \s)$,
\[
h(\lambda_0, \s^{p_\Cc})=p_\Cc\,h(\lambda, \s) \quad \text{and} \quad \int_{\SS_{\Cc, 0}}S_{p_\Cc}\p\,d\lambda_0=p_\Cc\int_{\SS_{\Cc}}\p\,d\lambda.
\]
Therefore a probability measure $\lambda \in \Mcc(\s, \SS_{\Cc})$ attains the supremum in \eqref{Eq:variational_principle}
if and only if the corresponding probability measure $\lambda_0 \in \Mcc(\s^{p_\Cc}, \SS_{\Cc, 0})$ attains the supremum
\eqref{Eq:Gibbs_ij}.
If we define $\m_\Cc:=\frac{1}{p_\Cc}\sum_{j=0}^{p_\Cc-1}\m_{\Cc, j}$,
then $\m_\Cc$ is the unique $\s$-invariant probability measure which attains the supremum of \eqref{Eq:variational_principle}.
Moreover, $\m_\Cc$ satisfies \eqref{Eq:Gibbs_i}.
This follows from \eqref{Eq:Gibbs-ij} and the H\"older continuity of $\p$.
Indeed,
for any $n \ge 1$, let $N$ be the integer with $Np_\Cc \le n < (N+1)p_\Cc$.
Taking
\[
[\o_0, \dots, \o_{(N+1)p_\Cc-1}] \subset [\o_0, \dots, \o_{n-1}] \subset [\o_0, \dots, \o_{Np_\Cc-1}],
\]
we have
\[
\m_\Cc[\o_0, \dots, \o_{n-1}] \le \frac{c_2}{p_\Cc} \cdot \exp\(p_\Cc(|\Pr_{\Cc}(\p, \s)|+\|\p\|_\infty)\)\cdot \exp\(-n \Pr_{\Cc}(\p, \s)+S_n \p(\o)\),
\]
for any $\o \in [\o_0, \dots, \o_{n-1}]$,
and
\[
\m_\Cc[\o_0, \dots, \o_{n-1}] \ge \frac{c_1}{p_\Cc} \cdot \exp\(-p_\Cc(|\Pr_{\Cc}(\p, \s)|+\|\p\|_\infty)-C'\)\cdot \exp\(-n \Pr_{\Cc}(\p, \s)+S_n \p(\o)\),
\]
for any $\o \in [\o_0, \dots, \o_{n-1}]$,
where 
$|S_{(N+1)p_\Cc}\p(\o')-S_n \p(\o)| \le C'+p_\Cc \|\p\|_\infty$ for any $\o' \in [\o_0, \dots, \o_{(N+1)p_\Cc-1}]$
and any $\o \in [\o_0, \dots, \o_{n-1}]$ for some constant $C' \ge 0$ independent of the cylinder sets and $\o$.
\qed

\subsection{Coding Patterson-Sullivan measures}

Let us take a H\"older continuous potential $\p:=\wh v\cdot \wt b$ on $\SS$, and let $\Cc_i$, $i=1, \dots, I$ be the maximal components of $\p$.
For each $i=1, \dots, I$, 
let $\m_i$ be a unique $\s$-invariant probability measure on $\SS_{\Cc_i}$ satisfying \eqref{Eq:Gibbs_i} in Proposition \ref{Prop:variational_principle}.
Let $\wt \m:=\sum_{i=1}^I \m_i$ on $\SS$.

\begin{lemma}\label{Lem:coding_PS}
Let $\m_{\wh d}$ be a Patterson-Sullivan measure for $\wh d$ on $\partial \G$.
There exists a constant $C>0$ such that 
\[
w_\ast \wt \m \le C \cdot \m_{\wh d} \otimes \m_{\wh d} \quad \text{on $\partial^2 \G$}.
\] 
Moreover, for any $T>0$, the measure 
\[
\sum_{\g \in \G}\g.(w_\ast \wt \m \otimes \Leb_{[0, T)})
\]
is Radon and absolutely continuous with respect to $\m_{\wh d} \otimes \m_{\wh d}\otimes dt$ on $\partial^2 \G \times \R$.
\end{lemma}

\begin{remark}
Note that $w_\ast \wt \m$ has the support in a compact subset of $\partial^2 \G$ while $\m_{\wh d}\otimes \m_{\wh d}$ has the full support $\partial^2 \G$.
We do not have the reverse inequality in the claim.
\end{remark}

\proof
Fixing a large enough $R>0$, we define shadows $O_o(x, R)$ on $\Cay(\G, S)$.
For any fixed $x, y \in \G$, let us consider any paths $\o \in \SS$ satisfying that for some $m, n \ge 0$,
\begin{equation}\label{Eq:codingPS}
\x_{-m}(\o) \in B_S(x, R) \quad \text{and} \quad \x_n(\o) \in B_S(y, R).
\end{equation}
Note that if some $\o$ in $[\o_{-m}, \dots, \o_n]$ satisfies \eqref{Eq:codingPS} with the indices $m, n\ge 0$, then in fact every $\o$ in the same cylinder set satisfies \eqref{Eq:codingPS}, and thus it is not ambiguous to say that a cylinder set $[\o_{-m}, \dots, \o_n]$ satisfies \eqref{Eq:codingPS}.
For every $\o$ in $[\o_{-m}, \dots, \o_{n-1}]$ with \eqref{Eq:codingPS}, we have
\[
\wh b_o(s_0\cdots s_{n-1}, w_\ast \o)=-\wh d(o, s_0\cdots s_{n-1})\pm C_\d
\]
since $\wh d$ is rough geodesic (see \eqref{Eq:Ancona} in Section \ref{Sec:strongly_hyperbolic}), and thus
\[
S_n\p(\o)=\wh v\cdot \wh b_o(s_0\cdots s_{n-1}, w_\ast \o)=-\wh v \cdot \wh d(o, s_0\cdots s_{n-1})\pm C_\d=-\wh v \cdot \wh d(o, y)\pm C_{R, \d},
\]
where in the last estimates we have used the triangle inequality,
and 
\[
S_m \p(\s^{-m}\o)=-\wh v \cdot \wh d(o, s_{-m}\cdots s_{-1}) \pm C_\d=-\wh v \cdot \wh d(x, o)\pm C_{R, \d},
\]
for $s_i=w(\o_i)$ for $i=-m, \dots, n-1$.
Since by Proposition \ref{Prop:shadow}, 
\[
\m_{\wh d}(O_o(x, R))=\exp\(-\wh v \cdot \wh d(o, x)\pm C\),
\]
we obtain
for all $\o \in [\o_{-m}, \dots, \o_{n-1}]$ with \eqref{Eq:codingPS},
\[
\m_{\wh d}(O_o(x, R))=\exp\(S_m \p(\s^{-m}\o) \pm C_{R, \d}'\) \quad \text{and} \quad \m_{\wh d}(O_o(y, R))=\exp\(S_n\p(\o)\pm C_{R, \d}'\).
\]
Note that for any $x, y \in \G$, we have
\begin{align*}
w_\ast \wt \m\Big(O_o(x, R) \times O_o(y, R)\Big)
&=\wt \m\Big(\o \in \SS \ : \ w_\ast(\o) \in O_o(x, R) \times O_o(y, R)\Big)\\
&\le \sum_{[\o_{-m}, \dots, \o_{n-1}]}\wt \m[\o_{-m}, \dots, \o_{n-1}],
\end{align*}
where the summation runs over all cylinder sets $[\o_{-m}, \dots, \o_{n-1}]$ satisfying \eqref{Eq:codingPS} (where $m, n$ also vary).
Since the measures $\m_i$ are $\s$-invariant and satisfy \eqref{Eq:Gibbs_i} in Proposition \ref{Prop:variational_principle} and $\Pr(\p, \s)=0$ by Lemma \ref{Lem:potential},
we have that for all $\o \in [\o_{-m}, \dots, \o_{n-1}]$ with \eqref{Eq:codingPS},
\begin{align*}
\wt \m[\o_{-m}, \dots, \o_{n-1}]	
\le C e^{S_{n+m}\p(\s^{-m}\o)} = C e^{S_m \p(\s^{-m}\o)}\cdot e^{S_n \p(\o)}
\le C'\m_{\wh d}(O_o(x, R))\cdot \m_{\wh d}(O_o(y, R)).
\end{align*}
Note that for each pair $x, y \in \G$,
the number of all cylinder sets $[\o_{-m}, \dots, \o_{n-1}]$ with \eqref{Eq:codingPS} is at most $|B_S(o, R)|^2$.
Hence for any $x, y \in \G$, we obtain
\[
w_\ast \wt \m\Big(O_o(x, R) \times O_o(y, R)\Big) \le C'|B_S(o, R)|^2\cdot \m_{\wh d}(O_o(x, R))\cdot \m_{\wh d}(O_o(y, R)).
\]
(Note that the left hand side may be $0$ for some $x$ and $y$.)
Comparing shadows with balls (Lemma \ref{Lem:shadows-balls}),  
we have constants $C, C'>0$ such that
for any $(\x, \y) \in \partial^2 \G$ and for any $r, s>0$,
\[
w_\ast \wt \m \Big(B_\rho(\x, r) \times B_\rho(\y, s)\Big) \le C \cdot \m_{\wh d}\otimes \m_{\wh d}\Big(B_\rho(\x, C'r) \times B_\rho(\y, C's)\Big).
\]
In the following discussion, we use a metric bi-Lipschitz to $\rho^\e$ for some $\e \in (0, 1)$ in $\partial \G$ (Section \ref{Sec:guage}) 
and denote it by the same symbol $\rho$.
Let us consider $\rho_\times((\x_1, \y_1), (\x_2, \y_2))=\max\{\rho(\x_1, \x_2), \rho(\y_1, \y_2)\}$ in $\partial^2 \G$.
For any Borel set $A$ in $\partial^2 \G$ and for any $\e>0$, we take a compact set $K$ and an open set $U$ in $\partial^2 \G$ such that
both $w_\ast \wt \m$ and $\m_{\wh d}\otimes \m_{\wh d}$ have their measures of $U \setminus K$ less than $\e$.
Since $K$ is compact, there exists an $\e_0>0$ such that every ball (relative to $\rho_\times$) of radius less than $\e_0$ centered in $K$ is included in $U$.
Let us take a countable family of balls $B_{\rho_\times}(x_i, r_i)$ of radius $r_i \le \e_0$ covering $K$ such that $B_{\rho_\times}(x_i, r_i/2)$ are disjoint.
Then we have that
\begin{align*}
w_\ast \wt \m\(K\) 
\le \sum_{i}w_\ast \wt \m(B_{\rho_\times}(x_i, r_i)) 
\le C \sum_{i}\m_{\wh d}\otimes \m_{\wh d}\(B_{\rho_\times}(x_i, C' r_i)\).
\end{align*}
Noting that $\m_{\wh d}\otimes \m_{\wh d}$ is doubling relative to $\rho_\times$ by Lemma \ref{Lem:doubling} (whose proof is adapted to the metric $\rho$) and $B_\rho(x_i, r_i/2)$ are disjoint, 
we obtain a constant $C'>0$ such that the last term is at most
\begin{align*}
C'\sum_{i} \m_{\wh d}\otimes \m_{\wh d}\Big(B_{\rho_\times}(x_i, r_i/2)\Big)
=C'\m_{\wh d}\otimes \m_{\wh d}\Big(\bigcup_i B_{\rho_\times}(x_i, r_i/2)\Big)
\le C'\m_{\wh d}\otimes \m_{\wh d}\(U\).
\end{align*}
Letting $\e\to 0$, we have $w_\ast \wt \m(A) \le C' \m_{\wh d}\otimes \m_{\wh d}(A)$ for any Borel set $A$ in $\partial^2 \G$.
Since the constant $C'$ is independent of $A$,
we conclude the first claim.

Let us show the second claim.
Note that for any $T>0$ and for every $\g \in \G$, the measure 
$\g.(w_\ast \wt \m \otimes \Leb_{[0, T)})$ is absolutely continuous with respect to $\m_{\wh d}\otimes \m_{\wh d}\otimes dt$,
and the finite measure $w_\ast \wt \m \otimes \Leb_{[0, T)}$ is supported on a compact set in $\partial^2 \G \times \R$ and the $(\G, \a)$-action of $\G$ on $\partial^2\G \times \R$ is properly discontinuous.
Therefore for any $T>0$ the measure 
$\sum_{\g \in \G}\g.(w_\ast \wt \m \otimes \Leb_{[0, T)})$
is Radon and absolutely continuous with respect to $\m_{\wh d}\otimes \m_{\wh d}\otimes dt$.
\qed

For the potential $\p_S=v_S \cdot \wt b_S$ as in Lemma \ref{Lem:potential_word},
letting $\Cc_i$, $i=1, \dots, I$ be the corresponding maximal components of $\p_S$,
we have a unique $\s$-invariant Borel probability measure $\m_i$ on $\SS_{\Cc_i}$ in Proposition \ref{Prop:variational_principle}.
(In fact, this specific case allows us to describe each $\m_i$ as the invariant measure for a Markov chain on $\Cc_i$ \cite[Section 4.2]{CalegariFujiwara2010}.) 
Let $\wt \m_S:=\sum_{i=1}^I \m_i$ on $\SS$.
The following lemma is proved as in Lemma \ref{Lem:coding_PS}; we omit the proof.

\begin{lemma}\label{Lem:coding_PS_word}
Let $\m_S$ be a Patterson-Sullivan measure associated with a word metric $d_S$ on $\partial \G$.
There exists a constant $C>0$ such that
\[
w_\ast \wt \m_S \le C \cdot \m_S \otimes \m_S \quad \text{on $\partial^2 \G$}.
\]
Moreover, for any $T>0$, the measure $\sum_{\g \in \G}\g.(w_\ast \m_S \otimes \Leb_{[0, T)})$ is Radon and absolutely continuous with respect to $\m_S \otimes \m_S \otimes dt$ on $\partial^2 \G \times \R$.
\end{lemma}

\section{Dimension and entropy}\label{Sec:entropy}

\subsection{Dimension}

For each $(\x_-, \x_+) \in \partial^2 \G$, let $\x$ be a geodesic in the Cayley graph $\Cay(\G, S)$ such that $\x_{-n} \to \x_-$ and $\x_{n} \to \x_+$ as $n \to \infty$ in $\G \cup \partial \G$, respectively.
For any hyperbolic metric $d$ in $\Dc_\G$,
we define
\[
\chi(\x_-, \x_+; d):=\liminf_{n \to \infty}\frac{1}{2n}d(\x_{-n}, \x_{n}).
\]
Note that $\chi(\x_-, \x_+)$ is independent of the choice of a geodesic $\x$ 
since any two geodesics with the same pair of extreme points are within a bounded distance in $\Cay(\G, S)$ and $d$ is quasi-isometric to a word metric.
For any $d \in \Dc_\G$, there exist constants $L_1, L_2>0$ such that
\[
L_1 \le \chi(\x_-, \x_+; d) \le L_2 \quad \text{for all $(\x_-, \x_+) \in \partial^2 \G$}.
\]

\begin{lemma}\label{Lem:Lyapunov}
Let $d \in \Dc_\G$, and
let $\L$ be a $\G$-invariant Radon measure on $\partial^2 \G$.
For $\L$-almost every $(\x_-,\x_+) \in \partial^2 \G$, and for any geodesic $\x$ in $\Cay(\G, S)$ with extreme points $(\x_-, \x_+)$,
the limit exists
\begin{equation}\label{Eq:chi}
\chi(\x_-, \x_+; d)=\lim_{n \to \infty}\frac{1}{2n}d(\x_{-n}, \x_{n}).
\end{equation}
Moreover, if $\L$ is ergodic, then $\chi(\x_-, \x_+; d)$ is constant $\L$-almost everywhere.
\end{lemma}

\proof
For $(\x_-, \x_+) \in \partial^2 \G$,
if the limit exists along a geodesic $\x$ in \eqref{Eq:chi}, then in fact the limit exists along any geodesics with the same pairs of extreme points.
Let us take a geodesic $\phi_{\x_-, \x_+}$ in $\Cay(\G, S)$ such that 
the map $(\x_-, \x_+) \mapsto \phi_{\x_-, \x_+}$ is measurable.
(Choosing a lexicographic ordering in the space of geodesics endowed with the topology of pointwise convergence, we assign the minimum $\phi_{\x_-, \x_+}$ in this ordering for each $(\x_-, \x_+)$; then this map is Borel measurable.)
Then, note that
the set
\[
B:=\Big\{(\x_-, \x_+) \in \partial^2 \G \ : \ \chi(\x_-, \x_+; d)=\lim_{n \to \infty}\frac{1}{2n}d(\x_{-n}, \x_{n}) \Big\}
\]
is Borel measurable and $\G$-invariant in $\partial^2 \G$.
Proposition \ref{Prop:coding_measure} (ii) implies that
for any $\G$-invariant Radon measure $\L$ on $\partial^2 \G$, there exists a $\s$-invariant probability measure $\lambda$ on $\SS$
such that 
\begin{equation}\label{Eq:equiv}
c_1 \L \otimes dt \le \sum_{\g \in \G}\g.\(w_\ast \lambda \otimes \Leb_{[0, T_0)}\) \le c_2 \L \otimes dt \quad \text{on $\partial^2 \G \times \R$}
\end{equation}
for some constants $c_1, c_2, T_0>0$.
The Kingman subadditive ergodic theorem implies that 
there exist the limits
\[
\lim_{n \to \infty}\frac{1}{n}d(\x_0(\o), \x_n(\o)) \quad \text{and} \quad \lim_{n \to \infty}\frac{1}{n}d(\x_{-n}(\o), \x_0(\o)) \quad \text{for $\lambda$-almost every $\o \in \SS$}.
\]
Since $d \in \Dc_\G$ is quasi-isometric to a word metric, there exists a constant $C_{d, S} \ge 0$ such that for all $\o \in \SS$ and for all $n \ge 1$,
\[
0\le d(\x_{-n}(\o), o) + d(o, \xi_n(\o))-d(\x_{-n}(\o), \x_n(\o))=2(\x_{-n}(\o)|\x_n(\o))_o \le C_{d, S}.
\]
Therefore the limit exists
\[
\chi(\x_-(\o), \x_+(\o); d)=\lim_{n \to \infty}\frac{1}{2n}d\,(\x_{-n}(\o), \x_n(\o)) \quad \text{for $\lambda$-almost every $\o \in \SS$},
\]
and thus $\lambda (w_\ast^{-1}B)=1$.
Since $B$ is $\G$-invariant, if we apply the set $(\partial^2 \G \setminus B) \times [0, 1]$ to the first inequality in \eqref{Eq:equiv},
then we obtain $\L(\partial^2 \G \setminus B)=0$.
Hence \eqref{Eq:chi} holds for $\L$-almost every $(\x_-, \x_+) \in \partial^2 \G$.
Moreover, since $\chi(\x_-, \x_+; d)$ is $\G$-invariant on $\partial^2 \G$, if $\L$ is ergodic with respect to the $\G$-action, then $\chi(\x_-, \x_+; d)$ has to be a constant function $\L$-almost everywhere on $\partial^2 \G$.
\qed

\begin{lemma}\label{Lem:Lyapunov_alpha}
If we define the function $\wt b$ on $\SS$ for the metric $\wh d$,
then for every $\s$-invariant Borel probability measure $\lambda$ on $\SS$, we have
\[
\int_{\SS}\wt b\,d\lambda=-\int_{\SS}\chi(\x_-(\o), \x_+(\o); \wh d\,)\,d\lambda.
\]
\end{lemma}
\proof
Since $\wh d$ is hyperbolic and quasi-isometric to a word metric, for any $\o \in \SS$ and any $n \ge 1$, 
we have for some constants $L, C\ge 0$,
\[
\lbr\x_{-n}(\o)| \x_{n}(\o)\rbr_o \le L(\x_{-n}(\o)| \x_{n}(\o))_o+C,
\]
(cf.\ \cite[Proposition 5.5 (1)]{BonkSchramm})
and thus
\[
S_n \wt b(\s^{-n}\o)+S_n \wt b(\o)=-\wh d(o, \x_{-n}(\o))-\wh d(o, \x_{n}(\o))\pm C=-\wh d(\x_{-n}(\o), \x_{n}(\o)) \pm C'.
\]
Therefore 
it holds that
\[
\lim_{n \to \infty}\frac{1}{2n}\(S_n \wt b(\s^{-n}\o)+S_n \wt b(\o)\)=-\chi(\x_{-}(\o), \x_{+}(\o); \wh d\,) \quad \text{for $\lambda$-almost every $\o \in \SS$}.
\]
Since $\wt b$ is bounded measurable on $\SS$, the Lebesgue dominated convergence theorem yields
\[
\int_\SS\wt b\,d\lambda=-\int_\SS\chi(\x_-(\o), \x_+(\o); \wh d\,)\,d\lambda,
\]
as required.
\qed

\subsection{Entropy}\label{Subsec:entropy}
Let $X$ be a compact topological space, and $f:X \to X$ be a continuous map.
We denote by $\Mcc(f, X)$ the space of all $f$-invariant Borel probability measures $\n$ on $X$.
For any $\n \in \Mcc(f, X)$, for any finite Borel partition $\Pc=\{P_1, \dots, P_k\}$ of $X$, i.e., each $P_i$ is a Borel set and $X=\bigsqcup_{i=1}^k P_i$,
let us define the entropy of $\Pc$ by
\[
H(\n, \Pc):=-\sum_{i=1}^k \n(P_i)\log\n(P_i).
\]
We define the entropy of $(X, \n, f)$ relative to $\Pc$ by
\[
h(f, \n, \Pc):=\lim_{n \to \infty}\frac{1}{n}H\Big(\n, \bigvee_{i=0}^{n-1}f^{-i}\Pc\Big)=\inf_{n \ge 1}\frac{1}{n}H\Big(\n, \bigvee_{i=0}^{n-1}f^{-i}\Pc\Big),
\]
where $\bigvee_{i=0}^{n-1}f^{-i}\Pc:=\Pc \vee f^{-1}\Pc \vee \cdots \vee f^{-(n-1)}\Pc$ is the partition consisting of all sets of the form
\[
P_{i_1}\cap f^{-1}P_{i_2}\cap \cdots \cap f^{-(n-1)}P_{i_n}, \quad i_1, \dots, i_n \in \{1, \dots, k\},
\]
and the limit exists by the subadditivity of the function $n \mapsto H(\n, \bigvee_{i=0}^{n-1}f^{-i}\Pc)$.
The {\it measure theoretical entropy} of $f$ for $(X, \n)$ is defined by
\[
h(f, \n):=\sup_{\Pc}h(f, \n, \Pc),
\]
where the supremum is taken over all finite Borel partitions $\Pc$.
For the suspension flow $\{\s_t\}_{t \in \R}$ on $\Sus(\SS, r)$, we consider the time one map $\s_1$ and consider $h(\s_1, \n)$ for a flow invariant probability measure $\n$ on $\Sus(\SS, r)$.

For any Radon measure $\L$ on $\partial^2 \G$,
let
\[
D_\rho(\x_-, \x_+; \L):=\liminf_{r \to 0}\frac{\log \L\(B_\rho(\x_-, r) \times B_\rho(\x_+, r)\)}{\log r} \qquad \text{for $(\x_-, \x_+) \in \partial^2 \G$}.
\]

\begin{lemma}\label{Lem:Local_dimension}
Let $d \in \Dc_\G$ and $\rho$ be the associated gauge in $\partial \G$.
Let $\L$ be a $\G$-invariant Radon measure on $\partial^2 \G$, and $\lambda$ be a $\s$-invariant Borel probability measure on $\SS$
such that for a constant $C>0$, we have
\[
w_\ast \lambda \le C\cdot \L \quad \text{on $\partial^2 \G$}.
\]
If $D_{\rho}(\x_-, \x_+; \L) \ge D$ for $\L$-almost every $(\x_-, \x_+)$ in $\partial^2 \G$, 
then we have that
\[
h(\s, \lambda) +\frac{D}{2} \int_{\SS}\wt b\,d\lambda \ge 0.
\]
\end{lemma}

\proof
We identify $\Sus\(\SS, r\)$ with $Z=\{(\o, t) \in \SS \times \R \ : \ 0\le t <r(\o)\}$.
Letting $W_N$ be the set of paths $(\o_0, \dots, \o_{N-1})$ of length $N$ in $\Ac$,
we decompose $Z$ into 
\[
\Big\{(\o, t) \in Z \ : \ \o \in [\o_0, \dots, \o_{N-1}], \ 0 \le t < r(\o)\Big\}
\]
for each cylinder set $[\o_0, \dots, \o_{N-1}]$ corresponding to the path in $W_N$ and 
define $\Pc$ the resulting partition of $Z$.
Thereby in $\Sus(\SS, r)$ we obtain the partition which we denote by the same symbol $\Pc$.
Let us fix $0<\e_0 < \min_{\o \in \SS} r(\o)$, and $\s_\ast:=\s_{\e_0}$ for the flow $\{\s_t\}_{t \in \R}$ on $\Sus(\SS, r)$.
We consider $\s_\ast ^n\Pc \vee \cdots \vee \s_\ast^{-n} \Pc$.

For any $[\o, t] \in \Sus(\SS, r)$ and for any integer $n \ge 1$, 
let us denote by $\Pc_{[-n, n]}([\o, t])$ the set in $\s_\ast ^n\Pc \vee \cdots \vee \s_\ast^{-n} \Pc$ containing $[\o, t]$.
Fix the corresponding point $(\o, t) \in Z$.
Then for $w_\ast(\o)=(\x_-(\o), \x_+(\o)) \in \partial^2 \G$, take a $C$-rough geodesic $\phi$ in $(\G, d)$ with $d(\phi(0), o) \le C$ such that $\phi(-n) \to \x_-(\o)$ and $\phi(n)\to \x_+(\o)$ as $n \to \infty$, respectively (Lemma \ref{Lem:BonkSchramm}).
Note that $\phi$ and $\x(\o)$ are within a bounded distance where $\x(\o)$ is a geodesic in $\Cay(\G, S)$ since $d$ and $d_S$ are quasi-isometric.
Fixing a large enough $R>0$, we have that
\[
[\o', t'] \in \Pc_{[-n, n]}([\o, t]) \implies w_\ast(\o') \in O_o\(\phi(t-\e_0 n), R\) \times O_o\(\phi(t+\e_0 n), R\).
\]
This follows 
since $\x(\o)$ and $\x(\o')$ coincides at the coordinates $i=-Nm', \dots, Nm$ where
$m:=\max\{k \ge 0 \ : \ S_k r(\o) \le t+\e_0 n\}$ and
$m':=\min\{k \ge 0 \ : \ S_k r(\s^{-kN}\o) \ge -t+\e_0 n\}$.

By the assumption
we have $w_\ast \lambda \le C\cdot \L$.
Let $\n_\lambda:=\frac{1}{\lambda \otimes \Leb(Z)}(\lambda \otimes \Leb)|_Z$.
Then we have that for $T_0:=\sup_{\o \in \SS} r(\o)$,
\begin{align*}
\n_\lambda\Big(\Pc_{[-n, n]}([\o, t])\Big)	&\le \n_\lambda\Big(w_\ast(\o') \in O_o\(\phi(t-\e_0 n), R\) \times O_o\(\phi(t+\e_0 n), R\), \ \ |t-t'| \le T_0\Big)\\
									&\le C\cdot \Lambda\(O_o\(\phi(t-\e_0 n), R\) \times O_o\(\phi(t+\e_0 n), R\)\)\cdot (2T_0).
\end{align*}
Comparing shadows with balls, we have that for all $n \ge 1$,
\[
O_o(\phi(t-\e_0 n), R) \times O_o(\phi(t+\e_0 n), R)\subset B_{\rho}(\x_-(\o), C_R e^{-\e_0 n}) \times B_{\rho}(\x_+(\o), C_R e^{-\e_0 n}),
\]
we obtain
\[
\liminf_{n \to \infty}-\frac{1}{n}\log \n_\lambda\Big(\Pc_{[-n, n]}([\o, t])\Big)
\ge \e_0 \cdot D_{\rho}(\x_-(\o), \x_+(\o); \Lambda).
\]
Noting that $\lambda$ is $\{\s_t\}_{t \in \R}$-invariant, by the Fatou lemma we obtain
\begin{align*}
&h(\s_\ast, \n_\lambda, \Pc) =\lim_{n \to \infty}\frac{1}{2n}H\Big(\n_\lambda, \bigvee_{i=0}^{2n+1} \s_\ast^{-i}\Pc\Big)
=\lim_{n \to \infty}\frac{1}{2n}H\Big(\n_\lambda,  \bigvee_{i=-n}^n \s_\ast^{-i}\Pc\Big)\\
&=\liminf_{n \to \infty}-\frac{1}{2n}\int_{\Sus(\SS, r)} \log \n_\lambda\(\Pc_{[-n, n]}([\o, t])\)\,d\n_\lambda
\ge \frac{\e_0}{2}\int_{\Sus(\SS, r)} D_{\rho}(\x_-(\o), \x_+(\o); \Lambda)\,d\n_\lambda.
\end{align*}
We have that $h(\s_\ast, \n_\lambda) \ge h(\s_\ast, \n_\lambda, \Dc)$ by definition of the measure theoretic entropy.
Since $w_\ast \lambda \le C\cdot \Lambda$,
if $D_{\rho}(\x_-, \x_+; \L) \ge D$ for $\L$-almost every $(\x_-, \x_+)$,
then $D_{\rho}(\x_-, \x_+; \L) \ge D$ for $w_\ast \lambda$-almost every $(\x_-, \x_+)$ in $\partial^2 \G$,
and thus we have $h(\s_\ast, \n_\lambda) \ge \e_0 D/2$.

Then we apply the Abramov formula and its consequence (e.g., \cite[Theorem 4.1.4 and Corollary 4.1.10]{FisherHasselblatt})
\[
h(\s^N, \lambda)=h(\s_1, \n_\lambda)\cdot \int_{\SS}r\,d\lambda \quad \text{and}
\quad h(\s_{\e_0}, \n_\lambda)=\e_0 \cdot h(\s_1, \n_\lambda)\qquad \text{for $\lambda \in \Mcc(\s^N, \SS)$}.
\]
Note that $\lambda$ is $\s$-invariant on $\SS$; we have $h(\s^N, \lambda)=N\cdot h(\s, \lambda)$, and 
$\int_\SS r\,d\lambda=N \int_\SS \wt \a\,d\lambda$.
Moreover, by Lemma \ref{Lem:alpha-Holder}, we have that $\wt \a=-\wt b+u\circ \s -u$ and thus 
$\int_\SS \wt \a\,d\lambda=-\int_\SS \wt b\,d\lambda$.
Therefore we obtain
\[
h(\s, \lambda) \ge -\frac{D}{2}\cdot \int_{\SS}\wt b\,d\lambda.
\]
as required.
\qed

\subsection{Invariant Radon measures of maximal dimensions}

\begin{lemma}\label{Lem:decompose}
Let $\lambda$ be a $\s$-invariant Borel probability measure on $\SS$.
For any component $\Cc$ of the underlying graph in the automatic structure $\Ac$, 
letting $\SS_\Cc$ be the set of bilateral paths all the time staying in $\Cc$,
we have that
\[
\lambda\Big(\SS\setminus\bigsqcup_{\Cc}\SS_\Cc\Big)=0.
\]
\end{lemma}
\proof
Let us denote by $\SS_{\to, \Cc}$ the set of bilateral paths eventually staying in $\Cc$ (and never leaving $\Cc$).
We note that
\[
\SS=\bigsqcup_\Cc \SS_{\to, \Cc},
\]
since there is no loop in the components graph (which is a directed graph obtained by identifying each component $\Cc$ in $\Ac$ with a point).
The set $\SS_{\to, \Cc}$ is $\s$-invariant and $\SS_{\to, \Cc}$ are disjoint for different $\Cc$.
For any integer $K$, if we define
\[
\SS_{K, \Cc}:=\Big\{\o \in \SS \ : \ \text{$\o_i$ is an edge in $\Cc$ for all $i \ge K$}\Big\},
\]
then we have $\SS_{\to, \Cc}=\bigcup_{K=-\infty}^\infty \SS_{K, \Cc}$ and $\SS_\Cc=\bigcap_{K=-\infty}^\infty \SS_{K, \Cc}$.
Note that $\SS_{K, \Cc} \subset \SS_{K+1, \Cc} \subset \s^{-1}\SS_{K, \Cc}$,
and $\SS_{\to, \Cc}\setminus \SS_\Cc=\bigsqcup_{K=-\infty}^\infty \(\SS_{K+1, \Cc}\setminus \SS_{K, \Cc}\)$.
Hence for any $\s$-invariant probability measure $\lambda$ on $\SS_{\to, \Cc}$,
we have that for each $K \in \Z$,
\[
\lambda\(\SS_{K+1, \Cc}\setminus \SS_{K, \Cc}\)\le \lambda\(\s^{-1}\SS_{K, \Cc}\)-\lambda\(\SS_{K, \Cc}\)=0,
\]
and thus $\lambda\(\SS_{\to, \Cc}\setminus \SS_\Cc\)=0$ for each $\Cc$.
Therefore every $\s$-invariant probability measure $\lambda$ on $\SS=\bigsqcup_\Cc \SS_{\to, \Cc}$ is supported on 
$\bigsqcup_\Cc \SS_\Cc$.
\qed

\begin{theorem}\label{Thm:uniqueMMD_main}
Let $\G$ be a non-elementary hyperbolic group and consider a strongly hyperbolic metric $\wh d \in \Dc_\G$ in $\G$ and the associated gauge $\wh \rho$ in $\partial \G$.
For any $\G$-invariant Radon measure $\Lambda$ on $\partial^2 \G$, 
if $\underline{\dim}_H(\Lambda, \wh \rho_\times)=\dim_H(\partial^2 \G, \wh \rho_\times)$, 
then $\Lambda$ is a constant multiple of $\L_{\wh d\,}$. 
\end{theorem}

\proof
Let us denote by $\wh v:=\gr(\wh d\,)$ the exponential volume growth rate relative to $\wh d$.
For any $\G$-invariant Radom measure $\Lambda$ on $\partial^2 \G$,
there exists a $\s$-invariant probability measure $\lambda$ on $\SS$ such that
\[
c_1 \Lambda \otimes dt \le \sum_{\g \in \G}\g.(w_\ast \lambda \otimes \Leb_{[0, T_0)}) \le c_2 \Lambda \otimes dt,
\]
for some constants $c_1, c_2, T_0>0$ by Proposition \ref{Prop:coding_measure} (ii).
Lemma \ref{Lem:decompose} implies that $\lambda$ is supported on $\bigsqcup_{\Cc}\SS_\Cc$.
Letting $\p=\wh v\cdot \wt b$ and 
$\lambda=\sum_\Cc a_\Cc\cdot \lambda_\Cc$,
where $a_\Cc \ge 0$ with $\sum_\Cc a_\Cc=1$, and $\lambda_\Cc$ is a $\s$-invariant probability measure supported on $\SS_\Cc$ for each component $\Cc$,
we have
\[
h(\s, \lambda)+\int_{\SS}\p\,d\lambda=\sum_{\Cc}a_\Cc\(h(\s, \lambda_\Cc)+\int_{\SS}\p\,d\lambda_\Cc\).
\]
The assumption $\underline \dim_H(\Lambda, \wh \rho_\times)=2\wh v$ implies that $D_{\wh \rho}(\x_-, \x_+; \L) \ge 2\wh v$ for $\L$-almost every $(\x_-, \x_+)$ in $\partial^2 \G$.
Hence Lemma \ref{Lem:Local_dimension} yields
\[
h(\s, \lambda)+\int_{\SS}\p\,d\lambda=h(\s, \lambda) +\wh v \int_\SS \wt b\,d\lambda \ge 0.
\]
The variational principle (Proposition \ref{Prop:variational_principle}) implies that for each $\Cc$,
\begin{equation}\label{Eq:vp_C}
\Pr_\Cc(\p, \s) \ge h(\s, \lambda_\Cc)+\int_{\SS}\p\,d\lambda_\Cc,
\end{equation}
and this implies that
\[
\Pr(\p, \s)=\max_\Cc \Pr_\Cc(\p, \s) \ge \sum_\Cc a_\Cc \cdot \Pr_\Cc(\p, \s) \ge h(\s, \lambda) +\int_{\SS}\p\,d\lambda \ge 0.
\]
For the potential $\p=\wh v\cdot \wt b$, we have $\Pr(\p, \s)=\max_{\Cc}\Pr_\Cc(\p, \s)=0$ by Lemma \ref{Lem:potential}.
Therefore we obtain
\[
\sum_\Cc a_\Cc \cdot \Pr_\Cc(\p, \s)=0 \quad \text{and} \quad \Pr_\Cc(\p, \s) \le 0 \quad \text{for all $\Cc$}.
\]
If $\Cc$ is not maximal for $\p$, then $\Pr_\Cc(\p, \s)<0$ and thus $a_\Cc=0$ for all component $\Cc$ which is not maximal for $\p$.
Hence $\lambda$ is supported on $\bigsqcup_{i \in I}\SS_{\Cc_i}$ where $\Cc_i$ for $i \in I$ are maximal components for $\p$.
If $a_{\Cc_i}>0$ for some $i \in I$, then for such an $i \in I$, we have that
\[
\Pr_{\Cc_i}(\p, \s) = h(\s, \lambda_{\Cc_i})+\int_{\SS}\p\,d\lambda_{\Cc_i}.
\]
Since for each $\Cc_i$, there exists a unique $\s$-invariant probability measure $\m_i$ which attains the equality in \eqref{Eq:vp_C} (Proposition \ref{Prop:variational_principle}), $\lambda_{\Cc_i}=\m_i$.
This shows that $\lambda=\sum_\Cc a_\Cc\cdot \lambda_\Cc$ is absolutely continuous with respect to $\wt \m =\sum_{i \in I}\m_i$,
and thus
$w_\ast \lambda$ is absolutely continuous with respect to $w_\ast \wt \m$.
By Lemma \ref{Lem:coding_PS}, we know that $w_\ast \wt \m$ is absolutely continuous with respect to $\L_{\wh d}$.
Therefore $\Lambda \otimes dt$ is absolutely continuous with respect to $\L_{\wh d}\otimes dt$, 
which shows that $\Lambda$ is absolutely continuous with respect to $\L_{\wh d}$.
Note that $d\Lambda/d\L_{\wh d}$ is locally integrable since $\L$ and $\L_{\wh d}$ are Radon and $\G$-invariant.
This shows that $d\Lambda/d\L_{\wh d}=c$ for a constant $c$ almost everywhere relative to $\L_{\wh d}$ since $\L_{\wh d}$ is ergodic with respect to the $\G$-action on $\partial^2 \G$ by Corollary \ref{Cor:double-ergodicity}.
Therefore we obtain $\Lambda=c\cdot \L_{\wh d}$ for some $c>0$ as desired.
\qed

\section{Mean distortion for word metrics}\label{Sec:mean_distortion}

\subsection{Word metrics}
Let $S$ be a finite set of generators $S$ with $S=S^{-1}$.
In this section, we focus on a word metric $d_S$ and denote by $\rho_S$ the associated gauge in $\partial \G$.
Let us consider the two-sided shift space $(\SS, \s)$ based on an automatic structure $(\Ac, w, S)$ with respect to $S$.

\begin{lemma}\label{Lem:Local_dimension_word}
Let $\L$ be a $\G$-invariant Radon measure on $\partial^2 \G$, and $\lambda$ be a $\s$-invariant Borel probability measure on $\SS$
such that for a constant $C>0$, we have
\[
w_\ast \lambda \le C\cdot \L \quad \text{on $\partial^2 \G$}.
\]
If $D_{\rho_S}(\x_-, \x_+; \L) \ge D$ for $\L$-almost every $(\x_-, \x_+)$ in $\partial^2 \G$, 
then we have that
\[
h(\s, \lambda)  \ge \frac{D}{2}.
\]
\end{lemma}

\proof
The proof runs similarly as in Lemma \ref{Lem:Local_dimension}.
For any integer $n \ge 1$, let $[\o_{-n}, \dots, \o_{n}]$ be any cylinder set of $\SS$.
Fix a large enough $R>0$ and we define shadows $O_o(x, R)$ on $\Cay(\G, S)$.
If $\o \in [\o_{-n}, \dots, \o_{n}]$, 
then
\[
w_\ast (\o) \in O_o(\x_{-n}(\o), R) \times O_o(\x_{n}(\o), R).
\]
Therefore for any $[\o_{-n}, \dots, \o_{n}]$, we have that
\begin{align*}
\lambda[\o_{-n}, \dots, \o_{n}]	&\le \lambda \circ w_\ast^{-1}\(O_o(\x_{-n}(\o), R) \times O_o(\x_{n}(\o), R)\)\\
								&=(w_\ast\lambda) \(O_o(\x_{-n}(\o), R) \times O_o(\x_{n}(\o), R)\).
\end{align*}
By the assumption that $w_\ast \lambda \le C \cdot \L$, the last term is at most
\[
C\cdot \L\(O_o(\x_{-n}(\o), R) \times O_o(\x_{n}(\o), R)\).
\]
Comparing shadows with balls (Lemma \ref{Lem:shadows-balls}), we have for any $n \ge 1$ and for any $\o \in \SS$,
\[
O_o(\x_{-n}(\o), R) \subset B_{\rho_S}(\x_-(\o), C_R e^{-n}) \quad \text{and} \quad O_o(\x_{n}(\o), R) \subset B_{\rho_S}(\x_+(\o), C_R e^{-n}),
\]
for some constant $C_R>0$ depending only on $R$, 
and for any $\o \in \SS$,
\begin{align*}
&\liminf_{n \to \infty}-\frac{1}{2n}\log \lambda[\o_{-n}, \dots, \o_{n}] \\
&\ge \liminf_{n \to \infty}-\frac{1}{2n} \log \L\(B_{\rho_S}(\x_-(\o), C_R e^{-n}) \times B_{\rho_S}(\x_+(\o), C_R e^{-n})\)
= D_{\rho_S}(\x_-(\o), \x_+(\o)).
\end{align*}
Let $\Pc:=\{\Pc_e\}_{e \in E}$ where $\Pc_e:=\{\o \in \SS \ : \ \o_0=e\}$ is the partition of $\SS$ according to the set of alphabets $E$.
By the definition of measure theoretic entropy, since $\lambda$ is $\s$-invariant, we have 
\begin{align*}
h(\s, \lambda, \Pc)	=\lim_{n \to \infty}\frac{1}{2n}H\Big(\s, \lambda, \bigvee_{i=0}^{2n}\s^{-i}\Pc\Big)
					=\lim_{n \to \infty}-\frac{1}{2n}\sum_{[\o_{-n}, \dots, \o_{n}]}\lambda[\o_{-n}, \dots, \o_{n}]\log \lambda[\o_{-n}, \dots, \o_{n}],
\end{align*}
and by the Fatou lemma, 
\begin{align*}
\liminf_{n \to \infty}-\frac{1}{2n}\sum_{[\o_{-n}, \dots, \o_{n}]}\lambda[\o_{-n}, \dots, \o_{n}]\log [\o_{-n}, \dots, \o_{n}]
\ge \int_{\SS}\liminf_{n \to \infty}-\frac{1}{2n}\log \lambda[\o_{-n}, \dots, \o_{n}]d\lambda.
\end{align*}
Hence we obtain
\begin{align*}
h(\s, \lambda)	&\ge \int_{\SS}D_{\rho_S}(\x_-(\o), \x_+(\o))\, d\lambda \ge \frac{D}{2}
\end{align*}
if $D_{\rho_S}(\x_-, \x_+) \ge D/2$ for $\L$-almost every $(\x_-, \x_+)$ in $\partial^2 \G$.
\qed

Let $\gr(S):=\gr(d_S)$ be the exponential volume growth rate with respect to $d_S$ and $\L_S$ be the corresponding Bowen-Margulis current for $d_S$.
Note that $\dim_H(\partial^2 \G, \rho_{S, \times})=2 \gr(S)$.

\begin{theorem}\label{Thm:uniqueMMD_word}
Let $\G$ be a non-elementary hyperbolic group, $d_S$ be a word metric on $\G$ and $\rho_S$ be the associated gauge in $\partial \G$.
For any $\G$-invariant Radon measure $\Lambda$ on $\partial^2 \G$, 
if $\underline{\dim}_H(\Lambda, \rho_{S, \times})=2\gr(S)$, 
then $\Lambda$ is a constant multiple of $\L_{S}$. 
\end{theorem}

\proof
The proof is analogous to Theorem \ref{Thm:uniqueMMD_main} and we indicate the place specific to word metrics.
Let $v_S:=\gr(S)$.
For any $\G$-invariant Radom measure $\Lambda$ on $\partial^2 \G$,
we take a $\s$-invariant probability measure $\lambda$ on $\SS$ in Proposition \ref{Prop:coding_measure}
and $\lambda$ is supported on 
$\bigsqcup_\Cc \SS_\Cc$ by Lemma \ref{Lem:decompose}. 
Let $\p_S:=v_S\cdot \wt b_S$ where $\wt b_S\equiv -1$.
The assumption $\underline \dim_H(\Lambda, \rho_{S, \times})=2v_S$ implies that $D_{\rho_S}(\x_-, \x_+; \L) \ge 2v_S$ for $\L$-almost every $(\x_-, \x_+)$ in $\partial^2 \G$.
Hence Lemma \ref{Lem:Local_dimension_word} implies that
$h(\s, \lambda) \ge v_S$,
and thus
\[
h(\s, \lambda)+\int_{\SS}\p_S\,d\lambda \ge 0.
\]
For the potential $\p_S=v_S\cdot \wt b_S$, we have $\Pr(\p_S, \s)=\max_{\Cc}\Pr_\Cc(\p_S, \s)=0$ by Lemma \ref{Lem:potential_word}.
We use the measure $\wt \m_S$ on $\SS$ in Lemma \ref{Lem:coding_PS_word}.
The rest follows similarly as in Theorem \ref{Thm:uniqueMMD_main}.
\qed

\begin{lemma}\label{Lem:distortion}
Let $S$ and $S^\star$ be finite symmetric sets of generators in $\G$.
There exists a constant $\t^\star$ such that for $\m_{S}$-almost every point $\x$ in $\partial \G$ and for any geodesic ray $\x_n$ in $\Cay(\G, S)$ converging to $\x$,
we have
\[
\lim_{n \to \infty}\frac{1}{n}d_{S^\star}(o, \x_n)=\t^\star,
\]
where $\m_S$ is a Patterson-Sullivan measure for $d_S$.
\end{lemma}

\proof
Let $\L_{S}$ be the Bowen-Margulis current on $\partial^2 \G$ associated with $d_S$.
Then we have that for $\L_S$-almost every $(\x_-, \x_+) \in \partial^2 \G$ and for any geodesic $\x$ in $\Cay(\G, S)$ with extreme points $(\x_-, \x_+)$, 
the limit exists 
\[
\chi^+(\x_-, \x_+; d_{S^\star}):=\lim_{n \to \infty}\frac{1}{n}d_{S^\star}(o, \x_n).
\]
Indeed, if this convergence holds for some geodesic ray toward $\x$, then in fact it holds for any geodesic ray toward $\x$ since any two geodesic rays converging to the same extreme point are eventually within bounded distance, and
this shows that the set where this limit exits is $\G$-invariant; the rest follows as in Lemma \ref{Lem:Lyapunov}.
Since $\L_S$ is ergodic with respect to $\G$-action on $\partial^2 \G$ by Corollary \ref{Cor:double-ergodicity}, 
this $\chi^+(\x_-, \x_+; d_{S^\star})$ is a constant function $\L_S$-almost everywhere on $\partial^2 \G$.
Note that $\L_S$ is equivalent to $\m_S \otimes \m_S$ and the limit does not depend on the choice of geodesic $\x$ converging to $\x_+$.
Letting the constant $\t^\star=\chi^+(\x_-, \x_+; d_S)$, we obtain the claim.
\qed

\subsection{Mean distortion}

Recall that for any pair of finite symmetric sets of generators $S$ and $S^\star$ in $\G$, 
we have defined
\[
\t(S^\star/ S):=\liminf_{n \to \infty}\frac{1}{n}\E_{\Unif_{S, n}}|x_n|_{S^\star},
\]
where $x_n$ has the uniform distribution $\Unif_{S, n}$ on $\bS_n=\{x \in \G \ : \ |x|_S=n\}$.
We shall show that the liminf is actually the limit, and the weak law of large number holds.

\begin{theorem}\label{Thm:distortion}
Let $\G$ be a non-elementary hyperbolic group.
For any pair of finite symmetric sets of generators $S$ and $S^\star$ in $\G$,
we have that
\[
\t(S^\star/ S)=\lim_{n \to \infty}\frac{1}{n}\E_{\Unif_{S, n}}|x_n|_{S^\star},
\]
and
for any $\e>0$,
\begin{equation*}
\lim_{n \to \infty} \Unif_{S, n}\left\{x \in \bS_n \ : \ \frac{|d_{S_\ast}(o, x)-n \t(S^\star/S)|}{n} > \e\right\}=0.
\end{equation*}
\end{theorem}

\proof
By Lemma \ref{Lem:distortion}, there exists a constant $\t^\star$ such that 
for $\m_S$-almost every $\x$ in $\partial \G$ and for any geodesic ray $\x_n$ converging to $\x$ in $\Cay(\G, S)$, we have
$d_{S^\star}(o, \x_n)/n \to \t^\star$ as $n \to \infty$.
For any $\e>0$, let
\[
A_n:=\left\{x \in \bS_n \ : \ \frac{|d_{S_\ast}(o, x)-n \t^\star|}{n} > \e \right\}.
\]
Fix a large enough $R>0$, and define shadows $O_o(x, R)$ in $\Cay(\G, S)$.
We have that $c_1 \le \m_S(O_o(x, R))\cdot |\bS_n| \le c_2$ for all $x \in \bS_n$ and for all $n \ge 0$
by Proposition \ref{Prop:shadow} and Lemma \ref{Lem:regular-growth} where we apply to $d_S$.
Note that shadows $O_o(x, R)$ for $|x|_S=n$ cover the boundary $\partial \G$, and
each $\x \in \partial \G$ is included at most $|B_S(o, 4R)|$ shadows $O_o(x, R)$ with $|x|_S=n$.
Therefore we have that
\begin{equation}\label{Eq:Thm:distortion}
\frac{|A_n|}{|\bS_n|} \le C\sum_{x \in A_n}\m_S\(O_o(x, R)\) \le C'\m_S\(\bigcup_{x \in A_n}O_o(x, R)\).
\end{equation}
If $\x \in \bigcup_{x \in A_n}O_o(x, R)$, 
then $\x_n \in B_S(x, 2R)$ for some $x \in A_n$, and thus $|d_{S^\ast}(o, \x_n)-n \t^\star| \ge \e n-2LR$ where $L$ is a Lipschitz constant: $d_{S^\star} \le Ld_S$.
Hence Lemma \ref{Lem:distortion} implies that the last term in \eqref{Eq:Thm:distortion} tends to $0$ as $n \to \infty$,
and we obtain
$\Unif_{S, n}(A_n) \to 0$ as $n \to \infty$.
This shows that
for any $\e>0$, for all large enough $n$, we have
$(\t^\star-\e)n \le |x_n|_{S^\star} \le (\t^\star+\e)n$ with probability at least $1-\e$,
and thus
\[
(1-\e)(\t^\star-\e)n \le \E_{\Unif_{S, n}}|x_n|_{S^\star} \le (\t^\star+\e)n+L\e n.
\]
Therefore we obtain
\[
\frac{1}{n}\E_{\Unif_{S, n}}|x_n|_{S^\star} \to \t^\star \quad \text{ as $n \to \infty$},
\]
and thus $\t^\star=\t(S^\star/S)$.
We conclude the claim.
\qed

\begin{lemma}\label{Lem:dimension_word}
Fix a word metric $d_{S^\star}$ in $\G$.
Let $\m_S$ be a Patterson-Sullivan measure for $d_S$.
Then
\[
\dim_H\(\m_S, \rho_{S^\star}\)=\frac{\gr(S)}{\t(S^\star/S)},
\]
where $\rho_{S^\star}$ is the gauge in $\partial \G$ relative to $d_{S^\star}$.
\end{lemma}

\proof
For any $\x$ in $\partial \G$, let us take a geodesic ray $\x^\star_n$ from $o$ converging to $\x$ in $\Cay(\G, S^\star)$.
Then $\x^\star$ is a $(L, C)$-quasi-geodesic ray in $\Cay(\G, S)$, and there exists a geodesic ray $\x_n$ from $o$ toward $\x$ in $\Cay(\G, S)$
such that $\x^\star_n$ and $\x_n$ are within bounded Hausdorff distance in $\Cay(\G, S)$.
Namely, there exists $D \ge 0$ such that for any $n$ we have $d_S(\x^\star_n, \x_{k_n}) \le D$ for some $k_n$.
Lemma \ref{Lem:distortion} implies that for $\m_S$-almost every $\x$ in $\partial \G$,
we have that
\[
\t(S^\ast/S) \le \liminf_{n \to \infty}\frac{d_{S^\star}(o, \x_{k_n})}{d_S(o, \x_{k_n})}  \le \limsup_{n \to\infty}\frac{d_{S^\star}(o, \x_{k_n})}{d_S(o, \x_{k_n})} \le \t(S^\star/S),
\]
and thus since $d_S(\x^\star_n, \x_{k_n}) \le D$ and $d_{S^\star}(\x^\star_n, \x_{k_n}) \le L\cdot D$ for all $n \ge 1$,
\[
\lim_{n \to \infty}\frac{d_{S^\star}(o, \x_n^\star)}{d_S(o, \x_n^\star)} = \t(S^\star/S) \quad \text{for $\m_S$-almost every $\x$}.
\]
Now we have $d_{S^\star}(o, \x_n^\star)=n$ and $\m_S\(O_o(x, R)\)$ is comparable to $\exp\(-\gr(S)|x|_S\)$ for all $x \in \G$ by Proposition \ref{Prop:shadow},
it holds that
\begin{equation}\label{Eq:dimension_word}
\lim_{n \to \infty}-\frac{1}{n}\log \m_S\(O_o(\x_n^\star, R)\)=\lim_{n \to \infty} \gr(S)\frac{d_S(o, \x_n^\star)}{d_{S^\star}(o, \x_n^\star)}=\frac{\gr(S)}{\t(S^\star/S)} \quad \text{for $\m_S$-almost every $\x$}.
\end{equation}
Comparing shadows with balls by Lemma \ref{Lem:shadows-balls} and by the Frostman-type lemma (Lemma \ref{Lem:Frostman}), we obtain $\dim_H(\m_S, \rho_{S^\star})=\gr(S)/\t(S^\star/S)$.
\qed

\begin{theorem}\label{Thm:distortion-similarity}
Let $\G$ be a non-elementary hyperbolic group.
For any pair of finite symmetric sets of generators $S$ and $S^\star$ in $\G$,
it holds that
\[
\t(S^\star/S) \ge \frac{\gr(S)}{\gr(S^\star)},
\]
and $\t(S^\star/S)=\gr(S)/\gr(S^\star)$ if and only if $d_{S^\star}$ and $d_S$ are roughly similar, i.e.,
there exist constants $\t>0$ and $D \ge 0$ such that
\[
|d_{S^\star}(x, y)-\t d_S(x, y)| \le D \quad \text{for all $x, y \in \G$}.
\]
\end{theorem}

\proof
Since 
$\dim_H(\m_S, \rho_{S^\star}) \le \gr(S^\star)$,
Lemma \ref{Lem:dimension_word} shows that we have that
$\t(S^\star/S) \ge \gr(S)/\gr(S^\star)$.
If there exist $\t>0$ and $D \ge 0$ such that $|d_{S^\star}(x, y)-\t d_S(x, y)| \le D$ for all $x, y \in \G$,
then $\t(S^\ast/S)=\t$ and $\gr(S)=\t \cdot \gr(S^\star)$
and thus the equality $\t(S^\star/S)=\gr(S)/\gr(S^\star)$ holds.

Let us assume that $\t(S^\star/S)=\gr(S)/\gr(S^\star)$.
Since $\L_S$ and $\m_S \otimes \m_S$ are equivalent on $\partial^2 \G$ and by \eqref{Eq:dimension_word} in Lemma \ref{Lem:dimension_word},
the Frostman-type lemma (Lemma \ref{Lem:Frostman}) implies that
\[
\dim_H(\L_S, \rho_{S^\star \times})=\frac{2\gr(S)}{\t(S^\star/S)}=2\gr(S^\star).
\]
Therefore Theorem \ref{Thm:uniqueMMD_word} implies that $\L_S=c \cdot \L_{S^\star}$ for some constant $c>0$.
This shows that $\m_S$ and $\m_{S^\star}$ are mutually absolutely continuous and thus $\m_S=f \m_{S^\star}$ for a density function $f \ge 0$ on $\partial \G$.
We shall show that $f$ is uniformly bounded from above and away from $0$.
Since $\L_S=c\cdot \L_{S^\star}$, by the definition of $\L_S$ and $\L_{S^\star}$,
we have that the ratio between
\[
f(\x)f(\y) \quad \text{and} \quad \frac{e^{2 \gr(S^\star)(\x|\y)_o^\star}}{e^{2 \gr(S)(\x|\y)_o}}
\]
are uniformly bounded from above and away from $0$ for all $\x \neq \y$ where $(\cdot | \cdot)_o$ and $(\cdot | \cdot)_o^\star$ are Gromov products with respect to $d_S$ and $d_{S^\star}$, respectively.
If the density $f$ is not bounded on $B_{\rho_S}(\x, \e)$, then fixing $\y \neq \x$ with $f(\y)>0$, we have that $f(\x')$ is arbitrary large for $\x' \in B_{\rho_S}(\x, \e)$ for all small enough $\e>0$ while $(\x'|\y)_o^\star$ and $(\x'|\y)_o$ are bounded for $\x' \in B_{\rho_S}(\x, \e)$; this is absurd.
Hence $f$ has to be uniformly bounded from above and the same argument shows that $f$ is uniformly bounded away from $0$.
Therefore $\m_S$ and $\m_{S^\star}$ are comparable, and for a fixed $R>0$ the ratio between
$\m_S(O_o(x, R))$ and $\m_{S^\star}(O_o(x, R))$ are uniformly bounded from above and away from $0$.
By Proposition \ref{Prop:shadow}, we have that
\[
\gr(S^\star)|x|_{S^\star}=\gr(S)|x|_S +O(1) \quad \text{for all $x \in \G$}.
\]
This shows that $d_{S^\star}$ and $d_S$ are roughly similar, as required.
\qed

\subsection*{Acknowledgements} 
The author would like to thank Professor Vadim Kaimanovich for sharing the questions and his insight which leads this work as well as valuable comments on an earlier version of this paper, and Professors J\'er\'emie Brieussel, Thibault Godin, Jayadev Athreya, Ilya Gekhtman, Hidetoshi Masai, Takefumi Kondo, Hiroyasu Izeki and Koji Fujiwara for helpful discussions on this topic, and Professor R\'emi Coulon for informing their works.
He also would like to thank the organizers of the workshop ``Equilibrium states for dynamical systems arising from geometry'' held at American Institute of Mathematics, where he had beneficial discussions with participants.
The author is supported by JSPS Grant-in-Aid for Young Scientists (B) JP17K14178
and JST, ACT-X Grant Number JPMJAX190J, Japan.

\bibliographystyle{alpha}
\bibliography{flow}

\end{document}